\documentclass{amsart}

\usepackage{a4wide}
\usepackage{amsmath,amscd}%
\usepackage{amsfonts}%
\usepackage{amssymb}%
\usepackage{graphicx}
\usepackage[usenames,dvipsnames]{color}
\usepackage{subfigure}
\usepackage{enumerate}
\usepackage[noprefix]{nomencl}  
\makenomenclature 
\usepackage[heads=vee]{diagrams}

\newcommand{\kdifform}[2]{#1^{(#2)}} 
\newcommand{\kdifformh}[2]{#1^{(#2)}_{h}} 
\newcommand{\kchain}[2]{\mathbf{#1}_{(#2)}} 
\newcommand{\kcochain}[2]{\mathbf{#1}^{(#2)}} 
\newcommand{\ccomplex}[1]{{#1}} 
\newcommand{\reconstruction}{\mathcal{I}} 
\newcommand{\reduction}{\mathcal{R}} 
\newcommand{\kformspace}[1]{\Lambda^{#1}} 
\newcommand{\kformspaceh}[1]{\Lambda_h^{#1}} 
\newcommand{\kformspacedomain}[2]{\Lambda^{#1}(#2)} 
\newcommand{\kchainspacedomain}[2]{C_{#1}(#2)} 
\newcommand{\kcochainspace}[1]{C^{#1}} 
\newcommand{\kcochainspacedomain}[2]{C^{#1}(#2)} 
\newcommand{\incidenceboundary}[2]{\mathsf{E}_{(#1,#2)}} 
\newcommand{\incidencederivative}[2]{\mathsf{E}^{(#1,#2)}} 
\newcommand{\projection}{\pi_{h}} 
\newcommand{\pullback}{\Phi^{\star}} 

\newcommand{\ederiv}{\mathrm{d}} 
\newcommand{\ud}{\mathrm{d}} 
\newcommand{\p}{\partial}
\newcommand{\dederiv}{\delta} 
\newcommand{\spacemap}{\rightarrow} 
 
\newcommand{\grad}{\mathrm{grad}\,} 
\newcommand{\curl}{\mathrm{curl}\,} 
\renewcommand{\div}{\mathrm{div}\,} 

\newcommand{\duality}[2]{\langle #1, #2\rangle} 

\renewcommand{\eqref}[1]{(\ref{#1})} 
\newcommand{\figref}[1]{Figure~\ref{#1}} 
\newcommand{\secref}[1]{Section~\ref{#1}} 
\newcommand{\defref}[1]{Definition~\ref{#1}} 
\newcommand{\lemmaref}[1]{Lemma~\ref{#1}} 

\theoremstyle{plain}

\newtheorem{definition}{Definition}
\newtheorem{example}{Example}

\newtheorem{lemma}{Lemma}

\numberwithin{equation}{section}

\newcommand{\define}{\mathrel{\mathop:}=}

\newcommand{\norm}[1]{\Vert #1\Vert}
\newcommand{\tr}{\mathrm{tr\;}}

\newcommand{\ve}{\varepsilon}

\graphicspath{{figures/}}

\begin{document}
\title[Mixed Mimetic Spectral Element Method]{Mixed Mimetic Spectral Element Method for Stokes Flow: A pointwise divergence-free solution}
\author{Jasper Kreeft}
\address[]
	{Delft University of Technology, Faculty of Aerospace Engineering, \newline%
	\indent  Kluyverweg 2, 2629 HT Delft, The Netherlands.}%
\email[]{J.J.kreeft@gmail.com, M.I.Gerritsma@TUDelft.nl}%
\author{Marc Gerritsma}
\thanks{Jasper Kreeft is funded by STW Grant 10113.}
\thanks{This paper is in final form and no version of it will be submitted for
publication elsewhere.}
\date{\today}
\subjclass{Primary 65M70, 76D07; Secondary 12Y05, 13P20} %
\keywords{Stokes problem, mixed finite elements, exact mass conservation, spectral elements, mimetic discretization}%

\begin{abstract}
In this paper we apply the recently developed mimetic discretization method to the mixed formulation of the Stokes problem in terms of vorticity, velocity and pressure. The mimetic discretization presented in this paper and in \cite{kreeftpalhagerritsma2011} is a higher-order method for curvilinear quadrilaterals and hexahedrals. Fundamental is the underlying structure of oriented geometric objects, the relation between these objects through the boundary operator and how this defines the exterior derivative, representing the grad, curl and div, through the generalized Stokes theorem. The mimetic method presented here uses the language of differential $k$-forms with $k$-cochains as their discrete counterpart, and the relations between them in terms of the mimetic operators: reduction, reconstruction and projection. The reconstruction consists of the recently developed mimetic spectral interpolation functions. The most important result of the mimetic framework is the commutation between differentiation at the continuous level with that on the finite dimensional and discrete level. As a result operators like gradient, curl and divergence are discretized exactly. For Stokes flow, this implies a pointwise divergence-free solution. This is confirmed using a set of test cases on both Cartesian and curvilinear meshes. It will be shown that the method converges optimally for all admissible boundary conditions.
\end{abstract}

\maketitle

\section{Introduction}\label{sec:introduction}
We consider Stokes flow, which models a viscous, incompressible fluid flow in which the inertial forces are negligible with respect to the viscous forces, i.e. when the Reynolds number is very small, $Re\ll1$. Since $Re=UL/\nu$, small Reynolds numbers appear when either considering extremely small length scales, when dealing with a very viscous liquid or when one treats slow flows.
Despite the simple appearance of Stokes flow model, there exists a large number of numerical methods to simulate Stokes flow. They all reduce to two classes of either circumventing the Ladyshenskaya-Babu\v{s}ka-Brezzi (LBB) stability condition or satisfying this condition, \cite{francahughes1988}. The first class can roughly be split into two subclasses, one is the group of stabilized methods, see e.g. \cite{bochevdohrmanngunzburger,hughes1986} and the references therein, the other the group of least-squares methods, see e.g. \cite{bochevgunzburger,jiang}.

The class that satisfies the LBB condition is the group of compatible methods. In compatible methods discrete vector spaces are constructed such that they satisfy the discrete LBB condition. Best known are the curl conforming N\'ed\'elec \cite{nedelec1980} and divergence conforming Raviart-Thomas \cite{raviartthomas1977} and Brezzi-Douglas-Marini \cite{brezzidouglasmarini1985} spaces. A subclass of compatible methods consists of {\em mimetic methods}. Mimetic methods do not solely search for appropriate vector spaces, but try to mimic structures and symmetries of the continuous problem, see \cite{bochevhyman2006,brezzibuffa2010,kreeftpalhagerritsma2011,mattiussi2000,perot2011,subramanian2006,tonti1}. As a consequence of this mimicking, mimetic methods automatically preserve structures of the continuous formulation.

\begin{figure}[htbp]
\centering
\includegraphics[width=0.7\textwidth]{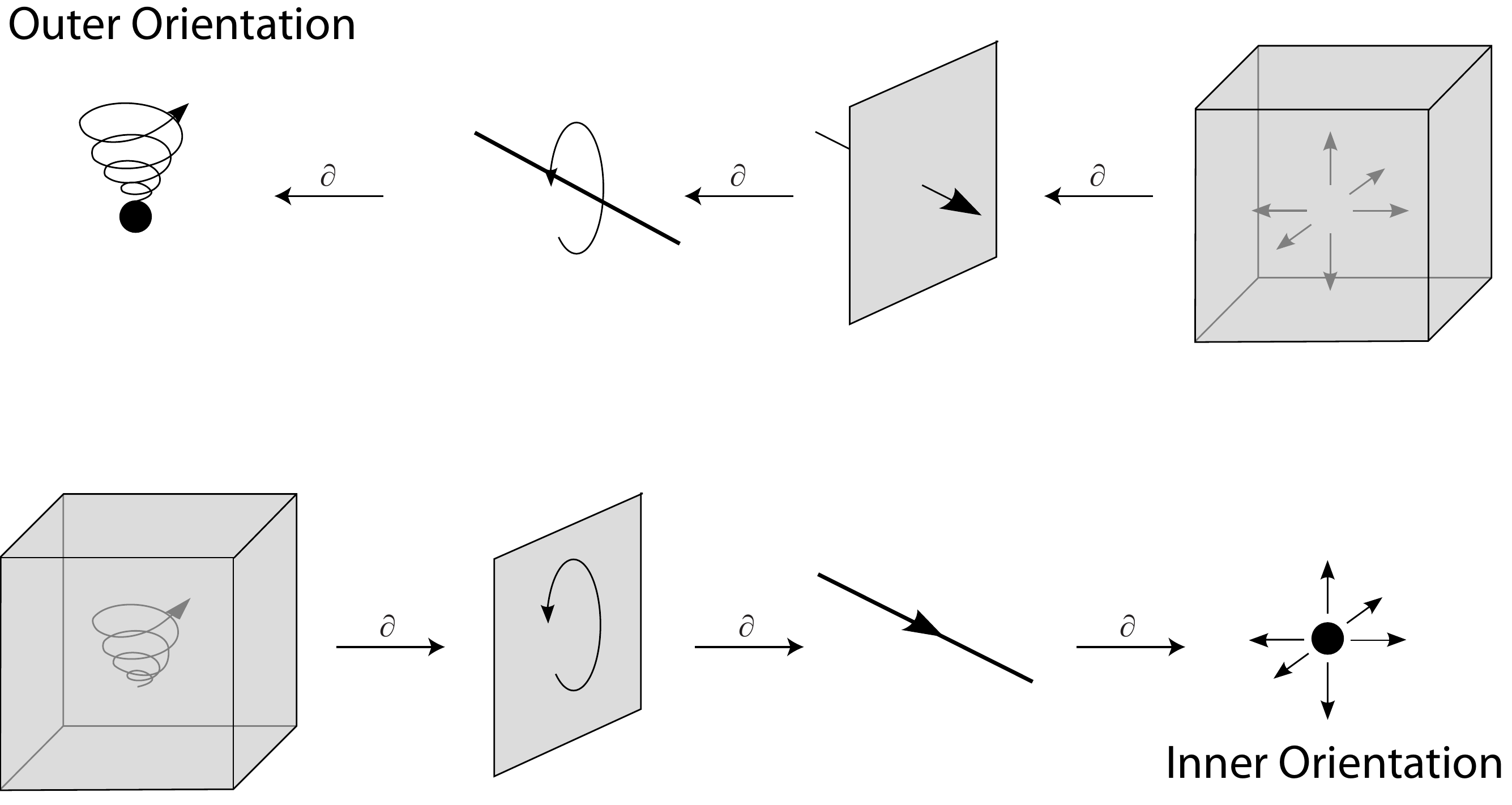}
\caption{The four geometric objects possible in $\mathbb{R}^3$, point, line, surface and volume, with outer- (above) and inner- (below) orientation. The boundary operator, $\partial$, maps $k$-dimensional objects to $(k-1)$-dimensional objects.}
\label{fig:manifoldswithorientation}
\end{figure}

At the heart of the mimetic method we present is the generalized Stokes theorem, which couples the exterior derivative to the boundary operator. In vector calculus this theorem is equivalent to the classical Newton-Leibnitz, Stokes circulation and Gauss divergence theorems. These well-known theorems relate the vector operators grad, curl and div to the restriction to the boundary of a manifold. Therefore, obeying geometry and orientation will result in satisfying exactly the mentioned theorems, and consequently performing the vector operators exactly in a finite dimensional setting. This is indeed what we are looking for and what our mimetic method has in common with finite volume methods, \cite{harlowwelch1965,subramanian2006}. In a three dimensional space we distinguish between four types of submanifolds, that is, points, lines, surfaces and volumes, and two types of orientation, namely, outer- and inner-orientation. The inner and outer orientations can be seen as generalizations of the concept of tangential and normal in vector calculus, respectively. This geometric structure will form the backbone of the mimetic method to be discussed in this paper. It will reappear throughout the paper in various guises. Examples of submanifolds in $\mathbb{R}^3$ are shown in \figref{fig:manifoldswithorientation} together with the action of the boundary operator.

By creating a quadrilateral or hexadedral mesh, we divide the physical domain in a large number of these geometric objects, and to each geometric object we associate a discrete unknown. This implies that these discrete unknowns are \emph{integral quantities}. Since the generalized Stokes theorem is an integral equation, it follows for example that taking a divergence in a volume is equivalent to taking the sum of the integral quantities associated to the surrounding surface elements, i.e. the fluxes. So using integral quantities as degrees of freedom to perform a vector operation like grad, curl or div, is equivalent to taking the sum of the degrees of freedom located at its boundary. These relations are of purely topological nature and can be seen as the horizontal connections between the geometric objects in \figref{fig:manifoldswithorientation}. The vertical connections -- not shown in \figref{fig:manifoldswithorientation} --, however describe the metric dependent parts, which are better known as the constitutive relations.


In this work we use the language of differential geometry to identify these structures, since it clearly identifies the metric and metric-free part of the PDEs. The latter has a discrete counterpart in the language of algebraic topology. In mimetic methods we employ commuting diagrams to indicate the strong analogy between differential geometry and algebraic topology. The most important commuting property employed in this work is the commutation between the projection operator and differentiation in terms of the exterior derivative. This means that also in finite dimensional spaces, operations like gradient, curl and divergence are performed exactly.
This implies, among others, and most importantly that incompressible Navier-Stokes and Stokes flow are guaranteed to be pointwise divergence-free, because the projection operator commutes with the divergence operator.

The similarities between differential geometry and algebraic topology in physical theories were first described by Tonti, \cite{tonti1}. A mimetic framework relating differential forms and cochains was initiated by Hyman and Scovel, \cite{HymanScovel1988}, and extended first by Bochev and Hyman, \cite{bochevhyman2006}, and later by Kreeft, Palha and Gerritsma \cite{kreeftpalhagerritsma2011}. A framework, closely related to the mentioned mimetic framework, is the finite element exterior calculus framework by Arnold, Falk and Winther \cite{arnoldfalkwinther2006,arnoldfalkwinther2010}. A more geometric approach is described in the work by Desbrun et al. \cite{desbrun2005c,desbrun2005}. An excellent introduction and motivation for the use of differential forms in the description of physics and the use in numerical modeling can be found in the `Japanese papers' by Bossavit, \cite{bossavit1998,bossavit9900}.

We make use of spectral element interpolation functions as basis functions. In the past nodal spectral elements were mostly used in combination with Galerkin (GSEM) \cite{bernardimayday,karniadakissherwin}, and least-squares formulations (LSSEM) \cite{pontaza2003,prootgerritsma2002}. The GSEM satisfies the LBB compatibility condition by lowering the polynomial degree of the pressure by two with respect to the velocity. This results in a method that is only weakly divergence-free, meaning that the divergence of the velocity field only convergence to zero with  mesh refinement. The LSSEM circumvents the LBB condition in order to be able to use equal order polynomials. The drawback of this method is the poor mass conservation property, \cite{kattelans2009,prootgerritsma2006}.

The present study uses mimetic spectral element interpolation or basis functions on curvilinear quadrilaterals and hexahedrals of arbitrary order as described in \cite{gerritsma2011,kreeftpalhagerritsma2011}. The mixed mimetic spectral element method (MMSEM) satisfies the LBB condition and gives a pointwise divergence-free solution for all mesh sizes. The mimetic spectral element interpolation functions are tensor product based interpolants. In every coordinate direction either a nodal or an edge interpolation function is used. By using tensor products, we are able to interpolate points, lines, surfaces, volumes, hyper-volumes and higher degree $n$-cube manifolds.

Although mimetic spectral elements are used to simulate Stokes flow and to derive numerical properties, alternative compatible/mimetic functions could be used in combination with the mimetic framework without much change, e.g., compatible B-splines, \cite{buffa2011,buffa2011b,evans2011}, and mimetic B-splines, \cite{back2011}.

This paper is organized as follows: first in \secref{sec:stokes} the Stokes problem in terms of vector calculus is given, with its relation to geometry and orientation. In \secref{sec:differentialgeometry} a brief summary of the most important concepts from differential geometry is given. \secref{sec:discretization} discusses the discretization of the Stokes model. It introduces the discrete structures of algebraic topology and a set of mimetic operators relating differential forms to cochains; the reduction operator, $\mathcal{R}$, the reconstruction operator, $\mathcal{I}$, and its composition, the bounded cochain projection, $\pi_h:=\mathcal{I}\circ\mathcal{R}$. As reconstruction functions the mimetic spectral element basis functions are used in this paper.
A mixed formulation for the Stokes problem is formulated in \secref{sec:mixedformulationstokes}. In \secref{sec:numericalresults} numerical results are discussed that show optimal convergence of all variables on curvilinear quadrilateral meshes. Secondly, the lid-driven cavity problem is shown on a square, cubic and triangle domain. The last test case is the flow around a cylinder moving with a constant velocity.

\section{Stokes problem in vector calculus}\label{sec:stokes}
Let $\Omega\subset\mathbb{R}^n$, $n=2,3$, be a bounded $n$-dimensional domain with boundary $\p\Omega$. On this domain we consider the Stokes problem, consisting of a momentum equation and the incompressibility constraint, resulting from the conservation of mass. The Stokes problem is given by
\begin{subequations}
\begin{align}
\nabla\cdot\sigma&=\vec{f}\quad\mathrm{on}\ \Omega,\\
\mathrm{div}\,\vec{u}&=0\quad\mathrm{on}\ \Omega,
\end{align}
\end{subequations}
where the stress tensor $\sigma$ is given by
\begin{equation}
\sigma=-\nu\nabla\vec{u}+p I,
\end{equation}
with $\vec{u}$ the velocity vector, $p$ the pressure, $\vec{f}$ the forcing term and $\nu$ the kinematic viscosity. In case of velocity boundary conditions the pressure is only determined up to a constant. So in a post processing step either the pressure in a point in $\Omega$ can be set, or a zero average pressure can be imposed; i.e.
\begin{equation}
\int_\Omega p\ud\Omega=0.
\end{equation}

For the method we like to present, we want to restrict ourselves to vector operations only. Therefore, instead of considering the divergence of a stress tensor, $\nabla\cdot(\nu\nabla\vec{u})$, we write this as $\nu\Delta\vec{u}$ by considering constant viscosity. Then the following vector identity is used for the vector Laplacian, $-\Delta\vec{u}=\mathrm{curl}\,\mathrm{curl}\,\vec{u}-\mathrm{grad}\,\mathrm{div}\,\vec{u}$. The vorticity-velocity-pressure formulation is obtained by introducing vorticity as auxiliary variable, $\vec{\omega}=\mathrm{curl}\,\vec{u}$. In terms of a system of first-order PDEs, the Stokes problem becomes
\begin{subequations}
\label{stokessinglevector}
\begin{align}
\vec{\omega}-\mathrm{curl}\,\vec{u}&=0\quad\mathrm{on}\ \Omega,\label{stokessinglevector1}\\
\mathrm{curl}\,\vec{\omega}+\mathrm{grad}\,p&=\vec{f}\quad\mathrm{on}\ \Omega,\label{stokessinglevector2}\\
\mathrm{div}\,\vec{u}&=0\quad\mathrm{on}\ \Omega.\label{stokessinglevector3}
\end{align}
\end{subequations}
Since these PDEs should hold on a certain physical domain, we will include geometry by means of integration. In that case we can relate every physical quantity to a geometric object. Starting with the incompressibility constraint \eqref{stokessinglevector3} we have due to Gauss' divergence theorem,
\[
\int_\mathcal{V}\div\vec{u}\,\ud\mathcal{V}=\int_{\p\mathcal{V}}\vec{u}\cdot\vec{n}\,\ud\mathcal{S}=0,
\]
and by means of Stokes' circulation theorem the relation \eqref{stokessinglevector1} can be written as
\[
\int_\mathcal{S}\vec{\omega}\times\vec{n}\,\ud\mathcal{S}=\int_\mathcal{S}\curl\vec{u}\times\vec{n}\,\ud\mathcal{S}=\int_{\p\mathcal{S}}\vec{u}\cdot\vec{t}\,\ud l.
\]
From the first relation it follows that $\div\vec{u}$ is associated to volumes. The association to a geometric object for velocity $\vec{u}$ is less clear. In fact it can be associated to two different types of geometric objects. In the incompressibility constraint velocity denotes a flux \emph{through} a surface that bounds a volume, while in the circulation relation velocity is defined \emph{along} a line that bounds the surface. We will call the velocity vector {\em through} a surface {\em outer-oriented} and the velocity {\em along} a line segment {\em inner-oriented}. Similarly, vorticity has two different representations, either as the rotation {\em in} a plane as Stokes' circulation theorem above suggests, or the Biot-Savart description of rotation {\em around} a line. In the former case $\vec{\omega}$ is {\em inner-oriented} whereas in the latter case $\vec{\omega}$ is {\em outer-oriented}, see also \figref{fig:manifoldswithorientation}. In fact both the velocity vector $\vec{u}$ and the vorticity vector $\vec{\omega}$ itself do not have a connection with geometry. Therefore, it are the terms $\vec{u}\cdot\vec{t}\,\ud l$, $\vec{u}\cdot\vec{n}\,\ud\mathcal{S}$, $\vec{\omega}\times\vec{n}\,\ud\mathcal{S}$ and $\vec{\omega}\times\vec{t}\,\ud l$ that are \emph{more useful variables} when considering Stokes problem on a physical domain.

\begin{figure}[htbp]
\centering
\includegraphics[width=0.7\textwidth]{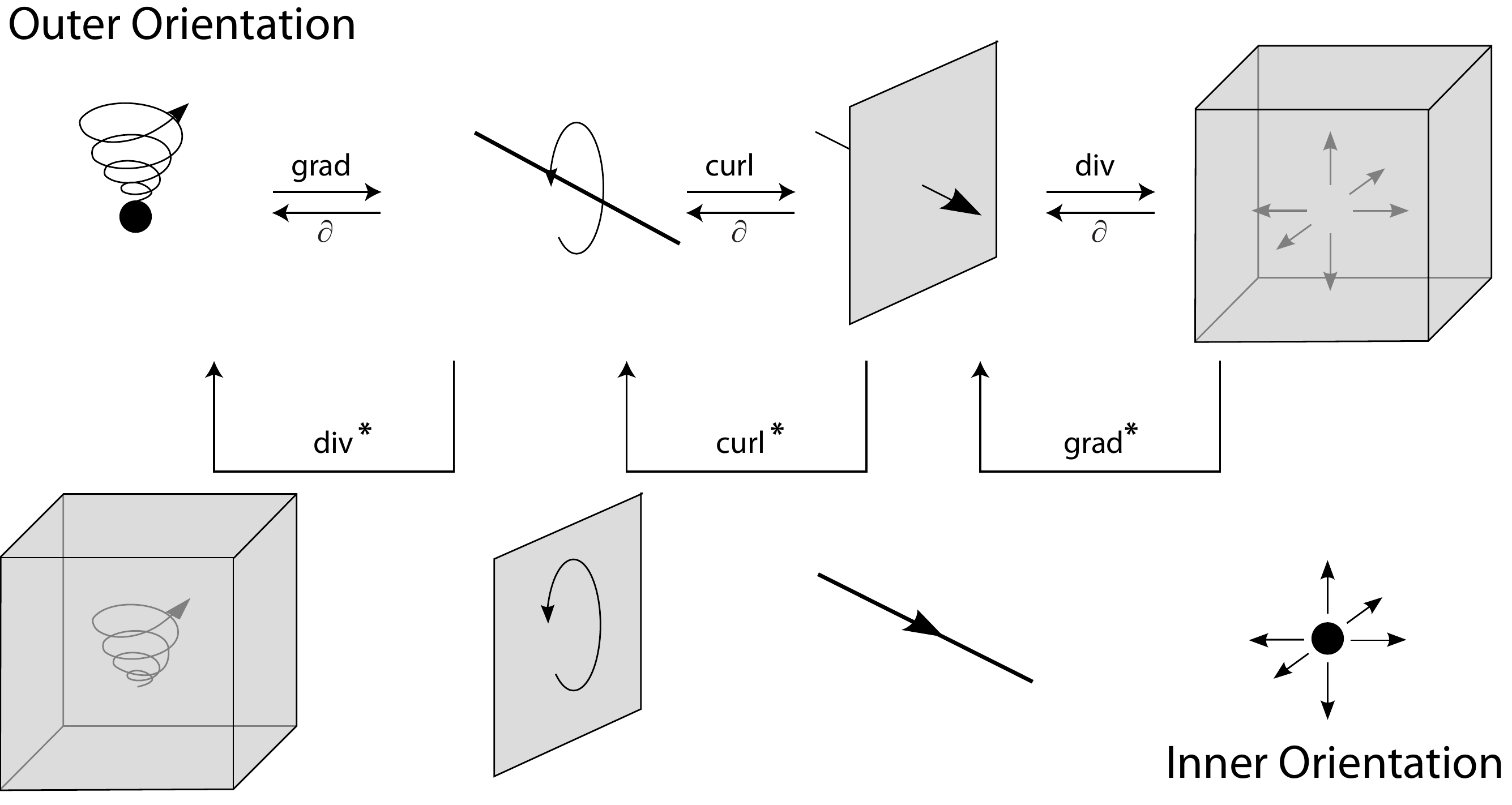}
\caption{Geometric interpretation of the action of the boundary operators, vector differential operators and their formal Hilbert adjoint operators.}
\label{fig:manifoldswithorientation2}
\end{figure}

The last equation to be considered is \eqref{stokessinglevector2}. This equation shows that classical Newton-Leibnitz, Stokes circulation and Gauss divergence theorems tell only half the story. From the perspective of the classical Newton-Leibnitz theorem, the gradient acting on the pressure relates line values to their corresponding end point, while the Stokes circulation theorem shows that the curl acting on the vorticity vector relates surface values to the line segment enclosing it. So how does this fit into one equation? In fact there exists two gradients, two curls and two divergence operators. One of each related to the mentioned theorems as explained above. The others are formal adjoint operators, so the second gradient is the adjoint of a divergence that is related to Gauss divergence theorem, the second curl is the adjoint of the curl related to Stokes circulation theorem and the second divergence is the adjoint of the gradient related to the classical Newton-Leibnitz theorem. Let grad, curl and div be the original differential operators associated to the mentioned theorems, then the formal Hilbert adjoint operators grad$^*$, curl$^*$ and div$^*$ are defined as,
\[
\big(\vec{a},\mathrm{grad}^*\,b\big):=\big(\div\vec{a},b\big),\quad \big(\vec{a},\mathrm{curl}^*\,\vec{b}\big):=\big(\curl\vec{a},\vec{b}\big),\quad \big(a,\mathrm{div}^*\,\vec{b}\big):=\big(\grad a,\vec{b}\big).
\]
Adjoint operators relate geometric operators in opposite direction. Where div relates a vector quantity associated to surfaces to a scalar quantity associated to a volume enclosed by these surfaces. Its adjoint operator, grad$^*$, relates a scalar quantity associated with a volume to a vector quantity associated with its surrounding surfaces. This is illustrated in \figref{fig:manifoldswithorientation2}. Following \figref{fig:manifoldswithorientation2}, the adjoint operators consist of three consecutive steps: First, switch to the other type of orientation (inner $\rightarrow$ outer or outer $\rightarrow$ inner), then take the derivative and finally switch the result back to its original orientation.

So \eqref{stokessinglevector2} could then either be associated to an inner-oriented line segment by rewriting it as
\[
\mathrm{curl}^*\,\vec{\omega}+\grad p=\vec{f},
\]
or be associated to an outer-oriented surface by rewriting it as
\[
\curl\vec{\omega}+\mathrm{grad}^* p=\vec{f}.
\]
Without geometric considerations we could never make a distinction between grad, curl and div and their associated Hilbert adjoints div$^*$, curl$^*$ and grad$^*$. Vector calculus does not make this distinction. 

Since our focus is on obtaining a pointwise divergence-free result, we decide to use a formulation associated to outer-oriented geometric objects. Then the Stokes problem becomes,
\begin{subequations}
\begin{align}
\vec{\omega}-\mathrm{curl}^*\,\vec{u}&=0,\\
\curl\vec{\omega}+\mathrm{grad}^*\,p&=\vec{f},\\
\div\vec{u}&=0,
\end{align}
\label{eq:system_vector_calculus}
\end{subequations}
where the first equation is associated to line segments, the second to surfaces and the third to volumes. In \cite{bochevgunzburger2009,bochevgunzburger} the same velocity-vorticity-pressure formulation is given in terms of grad, curl, div and grad$^*$, curl$^*$ and div$^*$.

For a valid equation, the mathematical objects should be the same; we can only add vectors with vectors and scalars with scalars, but not scalars with vectors. But now we add that equations also need the be \emph{geometrically compatible}. We can only add quantities associated with the same kind of geometry and with the same type of orientation. This lack of geometric notion in vector calculus is what motivates many to use the language of differential geometry instead, \cite{arnoldfalkwinther2006,arnoldfalkwinther2010,back2011,bochevhyman2006,bossavit1998,bossavit9900,buffa2011b,desbrun2005c,frankel,HymanScovel1988,kreefterrorestimate,kreeftpalhagerritsma2011,tonti1}. Other advantages of using differential geometry over vector calculus are that it possesses a clear distinction between variables associated with inner- and outer-orientation and it makes a clear distinction between topological and metric-dependent operations. All horizontal realtions in \figref{fig:manifoldswithorientation2} are topological. Any detour along geometric objects with the other type of orientation introduces metric in the equation. In differential geometry these structures are intrinsically embedded. It naturally leads to a discretization technique that can be seen as a hybrid between the finite volumes (topological part) and finite elements (metric part).

\section{Differential geometry}\label{sec:differentialgeometry}
This section presents the Stokes model in the language of differential forms. Differential geometry offers significant benefits in the construction of structure-preserving spatial discretizations. For example, the generalization of differentiation in terms of the exterior derivative encodes the gradient, curl and divergence operators from vector calculus and the codifferential represents the associated Hilbert adjoint operators grad$^*$, curl$^*$ and div$^*$. The generalized Stokes theorem encapsulates their corresponding integration theorems, respectively. The coordinate-free action of the exterior derivative and generalized Stokes theorem give rise to commuting properties with respect to mappings between different manifolds. These kind of commuting properties are essential for the structure preserving behavior of the mimetic method.

Only those concepts from differential geometry which play a role in the remainder of this paper will be explained. Much more can be found in \cite{abrahammarsdenratiu,flanders,frankel,kreeftpalhagerritsma2011}.

\subsection{Differential forms}
Let $\Lambda^k(\Omega)$ denote a space of \emph{differential $k$-forms} or \emph{$k$-forms}, on a sufficiently smooth bounded $n$-dimensional oriented manifold $\Omega\subset\mathbb{R}^n$, $\mathbf{x}:=(x^1,\hdots,x^n)$, with boundary $\p\Omega$. Every element $\kdifform{a}{k}\in\Lambda^k(\Omega)$ has a unique representation of the form,
\begin{equation}
\label{differentialform}
\kdifform{a}{k}=\sum_If_I(\mathbf{x})\ud x^{i_1}\wedge\ud x^{i_2}\wedge\cdots\wedge\ud x^{i_k},
\end{equation}
where $I=i_1,\hdots,i_k$, and $1\leq i_1<\hdots<i_k\leq n$ and where $f_I(\mathbf{x})$ are continuously differentiable scalar functions. Differential forms can be seen as quantities that live under the integral sign, \cite{flanders}, which were indicated in the previous section as the `more useful variables'. 
For $\Omega\subset\mathbb{R}^3$ with a Cartesian coordinate system $\mathbf{x}:=(x,y,z)$, the outer-oriented vorticity, velocity and pressure are
\begin{align*}
\kdifform{\omega}{1}&=\omega_1(\mathbf{x})\,\ud x+\omega_2(\mathbf{x})\,\ud y+\omega_3(\mathbf{x})\,\ud z,\\
\kdifform{u}{2}&=u(\mathbf{x})\,\ud y\!\wedge\!\ud z+v(\mathbf{x})\,\ud z\!\wedge\!\ud x+w(\mathbf{x})\,\ud x\!\wedge\!\ud y,\\
\kdifform{p}{3}&=p(\mathbf{x})\,\ud x\!\wedge\!\ud y\!\wedge\!\ud z.
\end{align*}
We see that $\omega$ is associated with line elements d$x$, d$y$ and d$z$. This is the outer-oriented representation in terms of Biot-Savart of rotation around a line segment. Similarly, velocity is associated with surface elements, $\ud y\!\wedge\!\ud z$, $\ud z\!\wedge\!\ud x$, $\ud x\!\wedge\!\ud y$, which is the outer-oriented representation of the velocity flux {\rm through} a surface. Finally, writing pressure as a volume form also corresponds to an outer-oriented representation.
 
Differential $k$-forms are naturally integrated over $k$-dimensional manifolds, i.e. for $\kdifform{a}{k}\in\Lambda^k(\Omega)$ and $\Omega_k\subset\mathbb{R}^n$, with $k=\mathrm{dim}(\Omega_k)$,
\begin{equation}
\label{integration}
\int_{\Omega_k}\kdifform{a}{k}\in\mathbb{R}\quad\Leftrightarrow\quad\langle \kdifform{a}{k},\Omega_k\rangle\in\mathbb{R},
\end{equation}
where $\langle\cdot,\cdot\rangle$ indicates a duality pairing\footnote{$\big(\cdot,\cdot\big)$ denotes metric dependent inner products, while $\langle\cdot,\cdot\rangle$ denotes metric-free duality pairing.} between the differential form and the geometry. This duality pairing is a metric-free operation, see \cite{frankel}. Note that the $n$-dimensional computational domain is indicated as $\Omega$, so without subscript. We would like to distinguish between $k$-forms that can be integrated over outer-oriented $k$-dimensional manifolds and $k$-forms that can be integrated over inner-oriented $k$-dimensional manifolds. To emphasis this difference, we sometimes write the space of the latter as $\tilde{\Lambda}^k(\Omega)$.

The wedge product, $\wedge$, of two differential forms $\kdifform{a}{k}$ and $\kdifform{b}{l}$ is a mapping: $\wedge:\Lambda^k(\Omega)\times\Lambda^l(\Omega)\rightarrow\Lambda^{k+l}(\Omega),\ k+l\leq n$. The wedge product is a skew-symmetric operator, i.e. $\kdifform{a}{k}\wedge \kdifform{b}{l}=(-1)^{kl}\kdifform{b}{l}\wedge \kdifform{a}{k}$. The pointwise inner-product of $k$-forms, $(\cdot,\cdot):\Lambda^k(\Omega)\times\Lambda^k(\Omega)\rightarrow\mathbb{R}$, is constructed using inner products of one-forms, that is based on the inner product on vector spaces, see \cite{flanders,frankel}.

The wedge product and inner product induce the Hodge-$\star$ operator, $\star:\Lambda^k(\Omega)\rightarrow\tilde{\Lambda}^{n-k}(\Omega)$, a metric operator that includes orientation, defined as
\begin{equation}
\kdifform{a}{k}\wedge\star \kdifform{b}{l}:=\big(\kdifform{a}{k},\kdifform{b}{k}\big)\kdifform{\sigma}{n},
\label{hodgestar}
\end{equation}
where $\kdifform{\sigma}{n}\in\Lambda^n(\Omega)$ is a unit volume form, $\kdifform{\sigma}{n}=\star1$. Let $(\ud x,\ud y,\ud z)$ be a basis in $\mathbb{R}^3$ for 1-forms associated with inner-oriented line segments. Then by applying the Hodge-$\star$ we retrieve a basis for 2-forms associated with outer-oriented surfaces,
\[
\star\ud x=\ud y\!\wedge\!\ud z,\quad \star\ud y=\ud z\!\wedge\!\ud x,\quad \star\ud z=\ud x\!\wedge\!\ud y.
\]
Therefore, the Hodge operator switches between inner- and outer-orientation. The Hodge-$\star$ operation can be interpreted as the vertical relations as given in \figref{fig:manifoldswithorientation2}, and coincides with a constitutive relation. The space of $k$-forms on $\Omega$ can be equipped with an $L^2$ inner product, $\big(\cdot,\cdot\big)_\Omega:\Lambda^k(\Omega)\times\Lambda^k(\Omega)\rightarrow\mathbb{R}$, given by,
\begin{equation}
\big(\kdifform{a}{k},\kdifform{b}{k}\big)_\Omega:=\int_\Omega\big(\kdifform{a}{k},\kdifform{b}{k}\big)\kdifform{\sigma}{n}=\int_\Omega \kdifform{a}{k}\wedge\star \kdifform{b}{k}.
\label{L2innerproduct}
\end{equation}
The differential forms live on manifolds and transform under the action of mappings. Let $\Phi:\Omega_{\rm ref}\rightarrow\Omega$ be a mapping between two manifolds. Then we can define the pullback operator, $\Phi^\star:\Lambda^k(\Omega)\rightarrow\Lambda^k(\Omega_{\rm ref})$, expressing the $k$-form on the reference manifold, $\Omega_{\rm ref}$. The mapping, $\Phi$, and the pullback, $\Phi^\star$, are related by
\begin{equation}
\int_{\Phi(\Omega_{\rm ref})}\kdifform{a}{k}=\int_{\Omega_{\rm ref}}\Phi^\star \kdifform{a}{k}\quad\Leftrightarrow\quad\langle \kdifform{a}{k},\Phi(\Omega_{\rm ref})\rangle=\langle\pullback \kdifform{a}{k},\Omega_{\rm ref}\rangle.
\end{equation}
A special case of the pullback operator is the trace operator. The trace of $k$-forms to the boundary, $\tr:\Lambda^k(\Omega)\rightarrow\Lambda^k(\p\Omega)$, is the pullback of the inclusion of the boundary of a manifold, $\p\Omega\hookrightarrow\Omega$, see \cite{kreeftpalhagerritsma2011}.

An important operator in differential geometry is the exterior derivative, $\ud:\Lambda^k(\Omega)\rightarrow\Lambda^{k+1}(\Omega)$. It represents the grad, curl and div (also rot in 2D) operators from vector calculus. It is induced by the \emph{generalized Stokes' theorem}, combining the classical Newton-Leibnitz, Stokes circulation and Gauss divergence theorems. Let $\Omega_{k+1}$ be a $(k+1)$-dimensional submanifold and $a^{(k)}\in\Lambda^k(\Omega)$, then
\begin{equation}\label{stokestheorem}
\int_{\Omega_{k+1}}\ud \kdifform{a}{k}=\int_{\p\Omega_{k+1}}\tr\kdifform{a}{k}\quad\Leftrightarrow\quad \langle\ud \kdifform{a}{k},\Omega_{k+1}\rangle = \langle\tr\kdifform{a}{k},\partial\Omega_{k+1}\rangle,
\end{equation}
where $\partial\Omega_{k+1}$ is a $k$-dimensional manifold being the boundary of $\Omega_{k+1}$. Due to the duality pairing in \eqref{stokestheorem}, the exterior derivative is the formal adjoint of the \emph{boundary operator} $\p:\Omega_{k+1}\rightarrow\Omega_k$ as indicated by the duality pairing, \eqref{integration}. The boundary operator \emph{defines} the exterior derivative.
The exterior derivative is independent of any metric or coordinate system. Applying the exterior derivative twice always leads to the null $(k+2)$-form, $\ud(\ud \kdifform{a}{k})=0^{(k+2)}$. On contractible domains the exterior derivative gives rise to an exact sequence, called \emph{de Rham complex} \cite{frankel}, and indicated by $(\Lambda,\ud)$,
\begin{equation}
\mathbb{R}\hookrightarrow\Lambda^0(\Omega)\stackrel{\ud}{\longrightarrow}\Lambda^1(\Omega)\stackrel{\ud}{\longrightarrow}\cdots\stackrel{\ud}{\longrightarrow}\Lambda^n(\Omega)\stackrel{\ud}{\longrightarrow}0.
\label{derhamcomplex}
\end{equation}
In vector calculus a similar sequence exists, where, from left to right for $\mathbb{R}^3$, the $\ud$'s denote the vector operators grad, curl and div. Both inner- and outer-oriented spaces of differential forms, $\Lambda^k(\Omega)$ and $\tilde{\Lambda}^k(\Omega)$, possess a de Rham sequence. The two are connected by the Hodge-$\star$ operator, and constitute a double de Rham complex,
\begin{equation}
\begin{matrix}
	\mathbb{R} \longrightarrow&\Lambda^{0}(\Omega)
	&\stackrel{\ederiv}{\longrightarrow}& \Lambda^{1}(\Omega)
	&\stackrel{\ederiv}{\longrightarrow}& \hdots \;
	&\stackrel{\ederiv}{\longrightarrow}\; &\Lambda^{n}(\Omega) \;
	&\stackrel{\ederiv}{\longrightarrow}\; &0  \\
	&\star\updownarrow & & \star\updownarrow &&  &
	&\star\updownarrow & &   \\
	0 \stackrel{\ederiv}{\longleftarrow}&{\tilde{\Lambda}}^{n}(\Omega)
	&\stackrel{\ederiv}{\longleftarrow}& {\tilde{\Lambda}}^{n-1}(\Omega)
	&\stackrel{\ederiv}{\longleftarrow}& \hdots \;
	&\stackrel{\ederiv}{\longleftarrow}\; &{\tilde{\Lambda}}^{0}(\Omega) \;
	&\stackrel{}{\longleftarrow}\; &\mathbb{R}.
\end{matrix} \label{double_deRham_complex}
\end{equation}
Observe the similarity between diagram \eqref{double_deRham_complex} and Figures \ref{fig:manifoldswithorientation} and \ref{fig:manifoldswithorientation2}, which is due to the fact that the exterior derivative is the adjoint of the boundary operator. The pullback operator and exterior derivative possess the following commuting property\footnote{Note that on the lefthandside of this equation we consider the pullback of a $(k+1)$-form, whereas on the right hand side the pullback of a $k$-form. We could have written this as $\Phi^\star_{k+1}\ud_k\kdifform{a}{k}=\ud_k\Phi^\star_k\kdifform{a}{k}$. In order to improve readibility and knowing that the meaning of these operators is clear from the context we do not explicitely denote this.},
\begin{equation}
\Phi^\star\ud \kdifform{a}{k}=\ud\Phi^\star \kdifform{a}{k},\quad \forall \kdifform{a}{k}\in\Lambda^k(\Omega),
\end{equation}
as illustrated in the following commuting diagram,
\[\begin{CD}
\Lambda^k(\Omega) @>\ud>> \Lambda^{k+1}(\Omega)\\
@VV\Phi^\star V @VV\Phi^\star V  \\
\Lambda^k(\Omega_{\rm ref}) @>\ud>> \Lambda^{k+1}(\Omega_{\rm ref}).
\end{CD}\]
The inner product, \eqref{L2innerproduct}, gives rise to the formal Hilbert adjoint of the exterior derivative, the {\emph codifferential operator}, $\ud^*:\Lambda^k(\Omega)\rightarrow\Lambda^{k-1}(\Omega)$, as $\big(\ud \kdifform{a}{k-1},\kdifform{b}{k}\big)_\Omega=\big(\kdifform{a}{k-1},\ud^* \kdifform{b}{k}\big)_\Omega$, which represents the grad$^*$, curl$^*$ and div$^*$ operators. Whereas the exterior derivative is a metric-free operator, the codifferential operator is metric-dependent, and given by $\ud^*=(-1)^{n(k+1)+1}\star\ud\star$, \cite{frankel,kreeftpalhagerritsma2011}. Here we see the three operations that were mentioned in the previous section and were illustrated in \figref{fig:manifoldswithorientation2}: Switch to the other type of orientation, $\star$, apply the derivative, d, and switch back to the original orientation, $\star$. In case of non-zero trace, and by combining \eqref{L2innerproduct} and \eqref{stokestheorem}, we get
\begin{equation}
\big(\kdifform{a}{k-1},\ud^*\kdifform{b}{k}\big)_\Omega=\big(\ud \kdifform{a}{k-1},\kdifform{b}{k})_\Omega-\int_{\p\Omega} \tr \kdifform{a}{k-1}\wedge \tr\star \kdifform{b}{k}.
\label{integrationbyparts}
\end{equation}
This is better known as integration by parts and is often used in finite element methods to avoid the codifferential.
Also for the codifferential, on contractible manifolds there exists an exact sequence,
\begin{equation}
0\stackrel{\ud^*}{\longleftarrow}\Lambda^0(\Omega)\stackrel{\ud^*}{\longleftarrow}\Lambda^1(\Omega)\stackrel{\ud^*}{\longleftarrow}\cdots\stackrel{\ud^*}{\longleftarrow}\Lambda^n(\Omega)\hookleftarrow\mathbb{R}.
\end{equation}
Finally, the Hodge-Laplace operator, $\Delta:\Lambda^k(\Omega)\rightarrow\Lambda^k(\Omega)$, is constructed as a composition of the exterior derivative and the codifferential operator,
\begin{equation}
\label{laplace}
-\Delta \kdifform{a}{k}:=(\ud^*\ud+\ud\ud^*)\kdifform{a}{k}.
\end{equation}

\subsection{Hilbert spaces}
Function spaces play an important role in the analysis of numerical methods. Of importance in this paper are the Hilbert spaces. On an oriented Riemannian manifold, we can define Hilbert spaces for differential forms. Let all $f_I(\mathbf{x})$ in \eqref{differentialform} be functions in $L^2(\Omega)$, then $\kdifform{a}{k}$ in \eqref{differentialform} is a $k$-form in the Hilbert space $L^2\Lambda^k(\Omega)$. The norm corresponding to the space $L^2\Lambda^k(\Omega)$ is $\Vert \kdifform{a}{k}\Vert_{L^2\Lambda^k}=\sqrt{(\kdifform{a}{k},\kdifform{a}{k})_\Omega}$. Although extension to higher Sobolev spaces are possible, we focus here on the Hilbert space corresponding to the exterior derivative. The Hilbert space $H\Lambda^k(\Omega)$ is defined by
\begin{equation}
H\Lambda^k(\Omega)=\{\kdifform{a}{k}\in L^2\Lambda^k(\Omega)\;|\;\ederiv \kdifform{a}{k}\in L^2\Lambda^{k+1}(\Omega)\},
\end{equation}
and the norm corresponding to $H\Lambda^k(\Omega)$ is defined as
\begin{equation}
\Vert \kdifform{a}{k}\Vert^2_{H\Lambda^k}:=\Vert \kdifform{a}{k}\Vert^2_{L^2\Lambda^k}+\Vert\ederiv \kdifform{a}{k}\Vert^2_{L^2\Lambda^{k+1}}.
\end{equation}
The Hilbert complex, $(H\Lambda,\ud)$, a special version of the de Rham complex, is the exact sequence of maps and spaces given by
\begin{equation}
\mathbb{R}\hookrightarrow H\Lambda^0(\Omega)\stackrel{\ederiv}{\longrightarrow} H\Lambda^1(\Omega)\stackrel{\ederiv}{\longrightarrow}\cdots\stackrel{\ederiv}{\longrightarrow} H\Lambda^n(\Omega)\stackrel{\ud}{\longrightarrow}0.
\end{equation}
In vector operations the Hilbert complex becomes for $\Omega\subset\mathbb{R}^3$,
\begin{equation}\label{3dcomplex}
H^1(\Omega)\stackrel{\rm grad}{\longrightarrow} H(\mathrm{curl},\Omega)\stackrel{\rm curl}{\longrightarrow}H(\mathrm{div},\Omega)\stackrel{\rm div}{\longrightarrow} L^2(\Omega),
\end{equation}
and for $\Omega\subset\mathbb{R}^2$, either
\begin{equation}\label{2dcomplexes}
H^1(\Omega)\stackrel{\rm grad}{\longrightarrow} H(\mathrm{rot},\Omega)\stackrel{\rm rot}{\longrightarrow}L^2(\Omega),\quad\mathrm{or}\quad
H^1(\Omega)\stackrel{\rm curl}{\longrightarrow} H(\mathrm{div},\Omega)\stackrel{\rm div}{\longrightarrow}L^2(\Omega).
\end{equation}
The two are related by the Hodge-$\star$ operator \eqref{hodgestar}, see \cite{palha2010},
\begin{equation}
\label{doublehilbertcomplex}
\begin{matrix}
H\Lambda^{0}(\Omega)\!\!&\!\!\stackrel{\ederiv}{\longrightarrow}\!\!&\!\! H\Lambda^{1}(\Omega)\!\!&\!\!\stackrel{\ederiv}{\longrightarrow}\!\!&\!\!L^2\Lambda^{2}(\Omega)\\
\star\updownarrow & & \star\updownarrow & &\star\updownarrow   \\
L^2\Lambda^{2}(\Omega)\!\!&\!\!\stackrel{\ederiv}{\longleftarrow}\!\!&\!\!H\Lambda^{1}(\Omega)\!\!&\!\!\stackrel{\ederiv}{\longleftarrow}\!\!&\!\!H\Lambda^{0}(\Omega)
\end{matrix}
\quad\Leftrightarrow\quad
\begin{matrix}
H^1(\Omega)\!\!&\!\!\stackrel{\mathrm{curl}}{\longrightarrow}\!\!&\!\!H(\mathrm{div},\Omega)\!\!&\!\!\stackrel{\mathrm{div}}{\longrightarrow}\!\!&\!\!L^2(\Omega)\\
\star\updownarrow & & \star\updownarrow & &\star\updownarrow   \\
L^2(\Omega)\!\!&\!\!\stackrel{\mathrm{rot}}{\longleftarrow}\!\!&\!\!H(\mathrm{rot},\Omega)\!\!&\!\!\stackrel{\mathrm{grad}}{\longleftarrow}\!\!&\!\!H^1(\Omega).
\end{matrix}
\end{equation}
A similar double Hilbert complex can be constructed in $\mathbb{R}^3$. Again note the similarities with these double Hilbert complexes and that of \eqref{double_deRham_complex} and geometric structure depicted in Figures \ref{fig:manifoldswithorientation} and \ref{fig:manifoldswithorientation2}.



\subsection{Stokes problem in terms of differential forms} 
The kind of form a variable has is directly related to the kind of manifold this variable can be integrated over. For example, from a physics point of view velocity is naturally integrated \emph{along} a line (streamline), a 1-manifold, indicating that velocity is a 1-form. However, looking at the incompressibility constraint, velocity in incompressible (Navier)-Stokes equations is usually associated to a flux \emph{through} a surface, indicating that velocity should be an $(n-1)$-form ($n=\mathrm{dim}(\Omega)$). The two are directly related by the Hodge duality, $u^{(n-1)}=\star \tilde{u}^{(1)}$, see\footnote{With $\tilde{\cdot}$ we indicate a variable contained in the lower complex of \eqref{double_deRham_complex}.} \eqref{double_deRham_complex}. The Hodge-$\star$ not only changes the corresponding type of integral domain, but also its orientation (along a line = inner, through a surface = outer).

Note that the Hodge-$\star$ is often combined with a constitutive relation. In that case the two variables have clearly a different meaning. In incompressible flow models, mass density plays the role of material property, so we actually have $(\rho u)^{(n-1)}=\star_\rho\tilde{u}^{(1)}$. Since mass density is assumed to be equal to one in incompressible (Navier)-Stokes, this difference is less obvious.

As for the velocity, also for pressure and vorticity there exists an inner and outer oriented version. The inner oriented variables are pressure, $\tilde{p}\in\Lambda^0(\Omega)$, associated to point values, vorticity and $\tilde{\omega}\in\Lambda^2(\Omega)$, associated to circulation in a surface. Alternatively, there exists the set of outer-oriented variables, being the pressure, $p\in\Lambda^n(\Omega)$, measured in a volume, and vorticity, $\omega\in\Lambda^{n-2}(\Omega)$, corresponding to circulation around a line (both in case of $\Omega\subset\mathbb{R}^3$).

Both sets, $(\tilde{p}^{(0)},\tilde{u}^{(1)},\tilde{\omega}^{(2)})$ and $(\omega^{(n-2)},u^{(n-1)},p^{(n)})$ are used in literature to derive mixed formulations and numerical schemes. For the former see \cite{abboud2011,bramble1994} and for the latter see \cite{bernardi2006,dubois2002}.

To obtain a pointwise divergence-free solution, the incompressibility constraint is leading, and therefore the set of outer-oriented variables are used in this paper, $(\kdifform{\omega}{n-2},\kdifform{u}{n-1},\kdifform{p}{n})$, with forcing term $\kdifform{f}{n-1}$. Then the Stokes problem in terms of differential forms becomes,
\begin{subequations}
\begin{align}
-\nu\Delta\kdifform{u}{n-1}+\ud^*\kdifform{p}{n}&=\kdifform{f}{n-1},\quad\mathrm{on}\ \Omega,\\
\ud \kdifform{u}{n-1}&=0,\quad\quad\quad\ \mathrm{on}\ \Omega,\label{mass}
\end{align}
\end{subequations}
where $\Delta$ is the Hodge-Laplacian defined by \eqref{laplace}. Vorticity is introduced as auxiliary variable to cast this system into a system of first-order equations. Substitution of \eqref{laplace} and the incompressibility constraint \eqref{mass}, gives the vorticity-velocity-pressure formulation in terms of differential forms,
\begin{subequations}
\label{stokeseq}
\begin{align}
\kdifform{\omega}{n-2}-\ud^* \kdifform{u}{n-1}&=0,\quad\quad\quad\ \mathrm{on}\ \Omega,\label{stokeseq1}\\
\nu\ud\kdifform{\omega}{n-2}+\ud^*\kdifform{p}{n}&=\kdifform{f}{n-1},\quad\mathrm{on}\ \Omega,\label{stokeseq2}\\
\ud \kdifform{u}{n-1}&=0,\quad\quad\quad\ \mathrm{on}\ \Omega.\label{stokeseq3}
\end{align}
\label{eq:system_diff_forms}
\end{subequations}
Note the resemblance of this system with \eqref{eq:system_vector_calculus}. Note also that whereas grad, curl and div are only defined in $\mathbb{R}^3$, \eqref{eq:system_diff_forms} is valid in $\mathbb{R}^n$ for all $n \geq 1$.

The actions of the exterior derivatives and codifferentials in this system are illustrated below for a two-dimensional domain.
\begin{example}[\textbf{2D Stokes problem}]\label{ex:2dstokes}
Let $\Omega\subset\mathbb{R}^2$, with Cartesian coordinates $\mathbf{x}:=(x,y)$, and let the two-dimensional de Rham complex be equivalent to the second complex in \eqref{2dcomplexes}. Then velocity is expressed as
\[
\kdifform{u}{1}=-v(\mathbf{x})\ud x+u(\mathbf{x})\ud y.
\]
Applying the exterior derivative gives us a 2-form, the divergence of velocity,
\[
\ud \kdifform{u}{1}=\left(\frac{\p u}{\p x}+\frac{\p v}{\p y}\right)\ud x\wedge\ud y.
\]
Vorticity is a 0-form, $\kdifform{\omega}{0}=\omega(\mathbf{x})\in\Lambda^0(\Omega)$, and the curl of vorticity gives,
\[
\ud\kdifform{\omega}{0}=\frac{\p\omega}{\p x}\ud x+\frac{\p\omega}{\p y}\ud y.
\]
The gradient of pressure, $\kdifform{p}{2}=p(\mathbf{x})\ud x\wedge\ud y\in\Lambda^2(\Omega)$, is the action of the codifferential,
\[
\ud^*\kdifform{p}{2}=-\frac{\p p}{\p y}\ud x+\frac{\p p}{\p x}\ud y.
\]
Then the momentum equation follows,
\[
-\left(-\frac{\p\omega}{\p x}+\frac{\p p}{\p y}\right)\ud x+\left(\frac{\p\omega}{\p y}+\frac{\p p}{\p x}\right)\ud y=-f_y(x,y)\ud x+f_x(x,y)\ud y.
\]
In a similar way the vorticity-velocity relation can be obtained.
\end{example}

\section{Discretization of Stokes problem}\label{sec:discretization}
The mimetic discretization of the Stokes problem consists of three parts. First, the discrete structure is described in terms of chains and cochains from algebraic topology, the discrete counterpart of differential geometry. This discrete structure mimics a lot of properties of differential geometry. Secondly, mimetic operators are introduced that relate the continuous formulation in terms of differential forms to the discrete representation based on cochains. Thirdly, mimetic spectral element basis functions are described which satisfy the structure defined in the algebraic topology and mimetic operators sections, Sections~\ref{sec:algeb_topol} and \ref{mimeticoperators}, respectively. The action of the exterior derivative, i.e., grad, curl and div, are shown, which leads among others to the divergence-free solution.

\subsection{Algebraic Topology}\label{sec:algeb_topol}
\begin{figure}[htbp]
\centering
\includegraphics[width=0.35\textwidth]{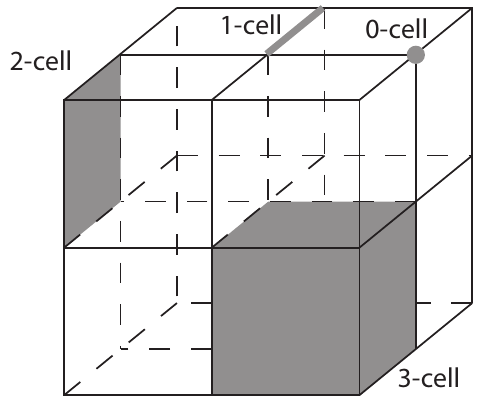}
\caption{Example of a 0-cell, a 1-cell, a 2-cell and a 3-cell in $\mathbb{R}^3$.}
\label{fig:kcells}
\end{figure}

In many numerical methods, especially in finite difference and almost all finite element methods, the discrete coefficients are point values, i.e. zero-dimensional sub-manifolds. In the structure of algebraic topology, the discrete unknowns represent values on $k$-dimensional submanifolds, ranging from points to $n$-dimensional volumes, so $0\leq k\leq n$. These $k$-dimensional submanifolds are called \emph{$k$-cells}, $\tau_{(k)}$. See \cite{hatcher,kreeftpalhagerritsma2011,munkres1984} how they are formally defined. The two most popular classes of $k$-cells in literature to describe the topology of a manifold are either in terms of \emph{simplices}, see for instance \cite{munkres1984,SingerThorpe,whitney}, or in terms of \emph{cubes}, see \cite{Massey2,spivak,tonti1} and \figref{fig:kcells} for an example of $k$-cubes in $\mathbb{R}^3$. From a topological point of view both descriptions are equivalent, see \cite{Dieudonne}. Despite this equivalence between simplicial complexes and cubical complexes, the reconstruction maps to be discussed in \secref{mimeticoperators} differ significantly. For mimetic methods based on simplices see \cite{arnoldfalkwinther2006,desbrun2005c,rapetti2009,subramanian2006}, whereas for mimetic methods based on cubes see \cite{arnoldboffifalk2005,HymanShashkovSteinberg2002,RobidouxSteinberg2011}.

Here we list the terminology to setup a homology theory in terms of $n$-cubes as given by \cite{Massey2}. Consider an oriented unit $k$-cube given by $I^k = I \times I \times \dots \times I$ ($k$ factors, $k\geq 0$), where $I=[0,1]$ is a one-dimensional closed interval. By definition $I^0$ is a space consisting of a single point. Then a {\em $k$-cube} in an $n$-dimensional manifold $\Omega$ is a continuous map $\tau_{(k)}:I^k\rightarrow\Omega,\ 0\leq k\leq n$.


All $k$-cells are oriented. This means that we define a \emph{default orientation}. The default orientation of the cell is implied by the orientation of the line segment $I$, which is defined positive in positive coordinate axis direction, and the map $\tau_{(k)}$. For outer-oriented cells, this for example also implies a positive way of going through a surface and rotating around a line. A $k$-cell with opposite orientation is said to have a negative orientation.

\begin{figure}[htbp]
\centering
\includegraphics[width=0.7\textwidth]{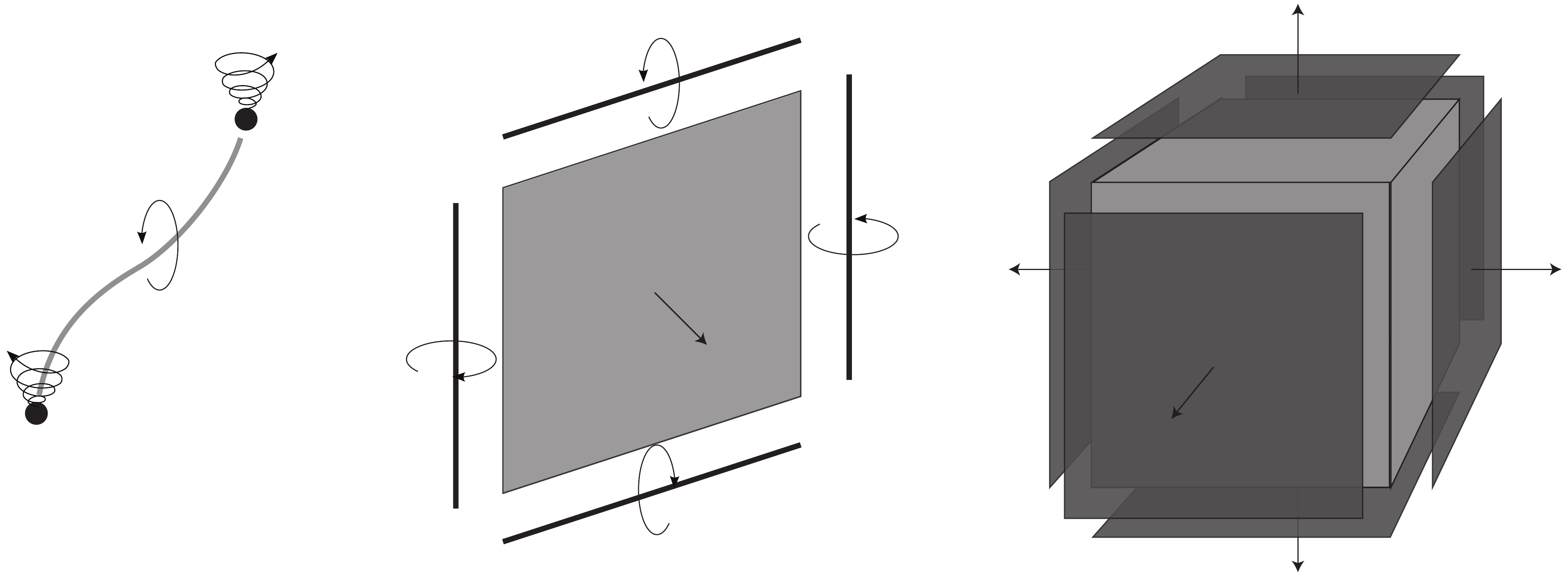}
\caption{Examples of faces of outer-oriented $k$-cells in $\mathbb{R}^3$.}
\label{fig:faces}
\end{figure}

The concept of orientation shown in Figures \ref{fig:manifoldswithorientation} and \ref{fig:manifoldswithorientation2} gives rise to the boundary operator, $\p$, that relates a $k$-cell to a set of surrounding $(k-1)$-cells, which has either the same or opposite orientation. Examples are given in \figref{fig:faces}, where the faces of the $k$-cells are shown in black.

This definition describes the boundary which we already encountered in \figref{fig:manifoldswithorientation} and \eqref{stokestheorem}. The boundary of a $k$-cell again consists of a set of $(k-1)$-cells, as illustrated in \figref{fig:faces}. From this we can define a \emph{cell complex}.
\begin{definition}[\textbf{Cell complex}]\cite{hatcher}\label{cellcomplex}
A cell complex, $\ccomplex{D}$, in a compact manifold $\Omega$ is a finite collection of cells such that:
\begin{enumerate}
\item The set of $n$-cells in $D$ covers the manifold $\Omega$.
\item Every face of a cell in $D$ is contained in $D$.
\item The intersection of any two $k$-cells, $\tau_{(k)}$ and $\sigma_{(k)}$ in $D$ either share a common $l$-cell, $l=0,\hdots,k-1$ in $D$, is empty, or $\tau_{(k)}=\sigma_{(k)}$.
\end{enumerate}
\end{definition}
\begin{figure}[htb]
\centering
\includegraphics[width=\textwidth]{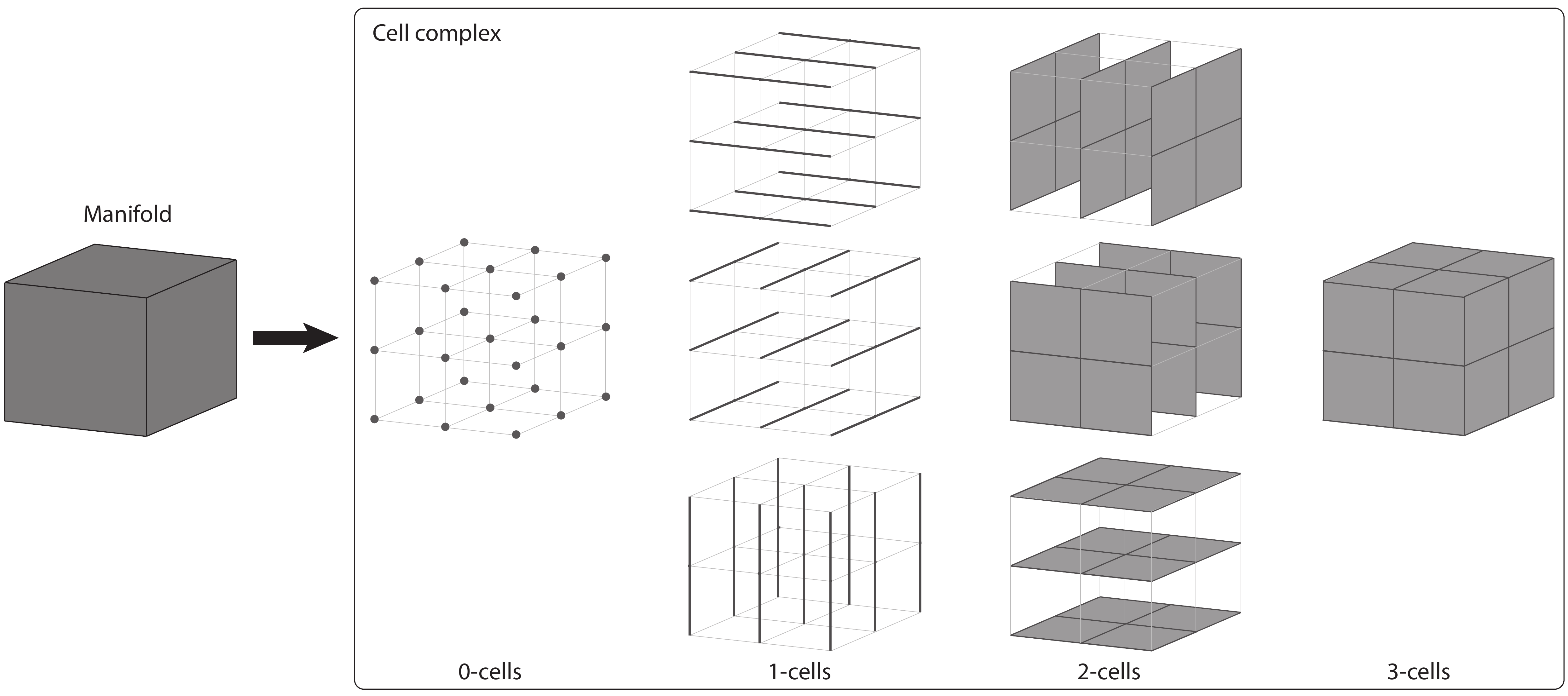}
\caption{Example of a cell complex. Left: a three dimensional compact manifold. Right: the $k$-cells that consistitute the cell complex.}
\label{fig:cellcomplex}
\end{figure}
We call a cell complex an \emph{oriented cell complex}, once we add to each $k$-cell a default orientation according to the definition of $k$-cubes. \figref{fig:cellcomplex} depicts an example of a cell complex in a compact manifold $\Omega\subset\mathbb{R}^3$. The ordered collection of all $k$-cells in $D$ generate a basis for the space of $k$-chains, $C_k(D)$. Then a $k$-chain, $\kchain{c}{k}\in C_k(D)$, is a formal linear combination of $k$-cells, $\tau_{(k),i}\in D$,
\begin{equation}
\kchain{c}{k}=\sum_ic_i\tau_{(k),i}.
\label{kchain}
\end{equation}
The $k$-cells, $\tau_{(k),i}$ form a basis for the $k$-chains. Once such a basis with default orientation has been chosen, any chain is completely determined by the coefficients $c_i$ which can be arranged in a column vector $\vec{c}=[c_1,c_2,\hdots]^T$.
In the description of geometry, we restrict ourselves to chains with coefficients in $\mathbb{Z}/3=\{ -1,0,1\}$. The meaning of these coefficients is : 1 if the cell is in the chain with the same orientation as its default orientation in the cell complex, -1 if the cell is in the chain with the opposite orientation to the default orientation in the cell complex and 0 if the cell is not part of the chain.

\begin{figure}[htb]
\centering
\includegraphics[width=0.7\textwidth]{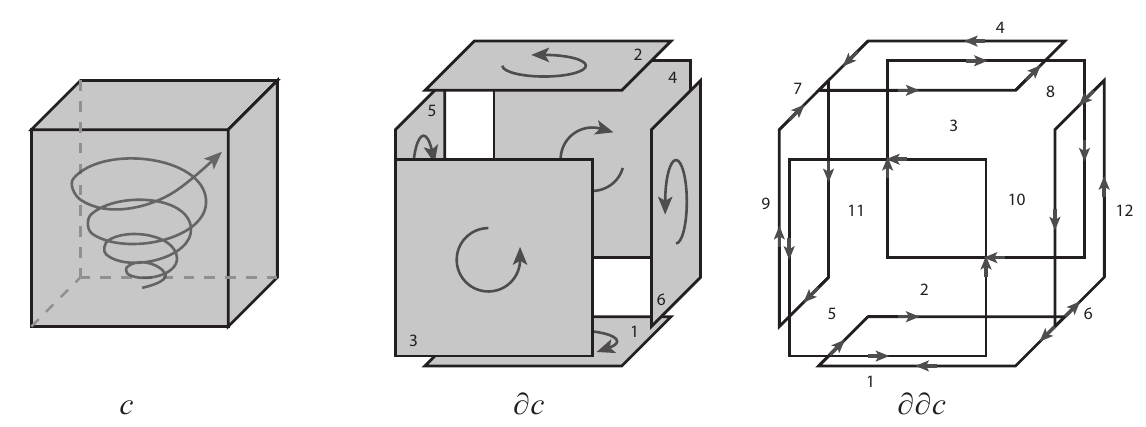}
\caption{The boundary of the boundary of a 3-cell is zero because all edges have opposite orientation.}
\label{fig:boundaryboundary}
\end{figure}

We can now extend the boundary operator applied to a $k$-cell to the boundary of a $k$-chain. The boundary operator acting on $k$-chains, $\partial:\kchainspacedomain{k}{D}\spacemap\kchainspacedomain{k-1}{D}$, is defined by \cite{hatcher,munkres1984},
\begin{equation}\label{algebraic::boundary_operator}
\partial \kchain{c}{k} = \partial \sum_{i}c^{i}\tau_{(k),i} := \sum_{i}c^{i} \partial \left ( \tau_{(k),i} \right ) \;.
\end{equation}
The boundary of a $k$-cell $\tau_{(k)}$ is a $(k-1)$-chain formed by the faces of $\tau_{(k)}$. The coefficients of this ($k-1$)-chain associated to each of the faces is given by the orientations.
\[
\partial \tau_{(k),i} = \sum_{j}e^{j}_{i}\tau_{(k-1),j} \;,
\]
with
\[
\left\{
\begin{array}{l}
e^{j}_{i} = 1, \text{ if the orientation of } \tau_{(k-1),j} \text{ equals the default orientation,} \\
e^{j}_{i} = -1, \text{ if the orientation of } \tau_{(k-1),j} \text{ is opposite to the default orientation,} \\
e^{j}_{i} = 0, \text{ if }  \tau_{(k-1),j}\text{ is not a face of } \tau_{(k),i}\;. \\
\end{array}
\right.
\]
The boundary of a 0-cell is empty. In case all $k$-cells in the chain $\kchain{c}{k}$ have positive orientation, so $c^i=1$, then
\begin{equation}
\partial \kchain{c}{k} = \sum_{i}\sum_je^{j}_{i}\tau_{(k-1),j}. \label{eq::algTop_boundary}
\end{equation}
Recalling that the space of $k$-chains is a linear vector space it follows that the boundary operator can be represented as a matrix acting on the column vector $\vec{c}$ of the $k$-chain. The coefficients $e^{j}_{i}$ are the coefficients of an incidence matrix $\incidenceboundary{k-1}{k}$ that represents the boundary operator. Like the exterior derivative, applying the boundary operator twice on a $k$-chain gives the null $(k-2)$-chain, $\p\p\kchain{c}{k}=\kchain{0}{k-2}$ for all $\kchain{c}{k}\in C_k(D)$, see \figref{fig:boundaryboundary}.
This was expected, since the exterior derivative and boundary operator are related according to the generalized Stokes theorem, \eqref{stokestheorem}. This property is reflected in the incidence matrices, since they are matrix representations of the topological boundary operators. Therefore $\incidenceboundary{k-2}{k-1}\incidenceboundary{k-1}{k}=0$, where for \figref{fig:boundaryboundary} we have
\[
\incidenceboundary{1}{2}=
\tiny{
\left[
\begin{array}{rrrrrr}
    -1  &   0  &   1  &   0  &   0  &   0 \\
     1  &   0  &   0  &  -1  &   0  &   0 \\
     0  &   1  &  -1  &   0  &   0  &   0 \\
     0  &  -1  &   0  &   1  &   0  &   0 \\
     1  &   0  &   0  &   0  &  -1  &   0 \\
    -1  &   0  &   0  &   0  &   0  &   1 \\
     0  &  -1  &   0  &   0  &   1  &   0 \\
     0  &   1  &   0  &   0  &   0  &  -1 \\
     0  &   0  &  -1  &   0  &   1  &   0 \\
     0  &   0  &   1  &   0  &   0  &  -1 \\
     0  &   0  &   0  &   1  &  -1  &   0 \\
     0  &   0  &   0  &  -1  &   0  &   1 \\
\end{array}
\right]
}, \quad\quad
\incidenceboundary{2}{3} =
\tiny{
\left[
\begin{array}{r}
-1 \\ 1 \\ -1 \\ 1 \\ -1 \\ 1 \\
\end{array}
\right]
}.
\]
The set of $k$-chains and boundary operators gives rise to an exact sequence, the chain complex $(C_k(D),\p)$,
\begin{equation}
\label{chaincomplex}
\begin{CD}
\cdots @<\p<< C_{k-1}(D) @<\p<< C_k(D) @<\p<< C_{k+1}(D) @<\p<< \cdots.
\end{CD}
\end{equation}
This sequence is the algebraic equivalent of \figref{fig:manifoldswithorientation}. Dual to the space of $k$-chains, $C_k(D)$, is the space of \emph{$k$-cochains}, $C^k(D)$, defined as the set of all linear functionals, $\kcochain{c}{k}:C_k(D)\rightarrow\mathbb{R}$. The duality is expressed using the duality pairing $\langle\kcochain{c}{k},\kchain{c}{k}\rangle:=\kcochain{c}{k}(\kchain{c}{k})$. Note the resemblance between this duality pairing and the integration of differential forms, see \eqref{integration}.

\begin{figure}[htb]
\centering
\includegraphics[width=0.8\textwidth]{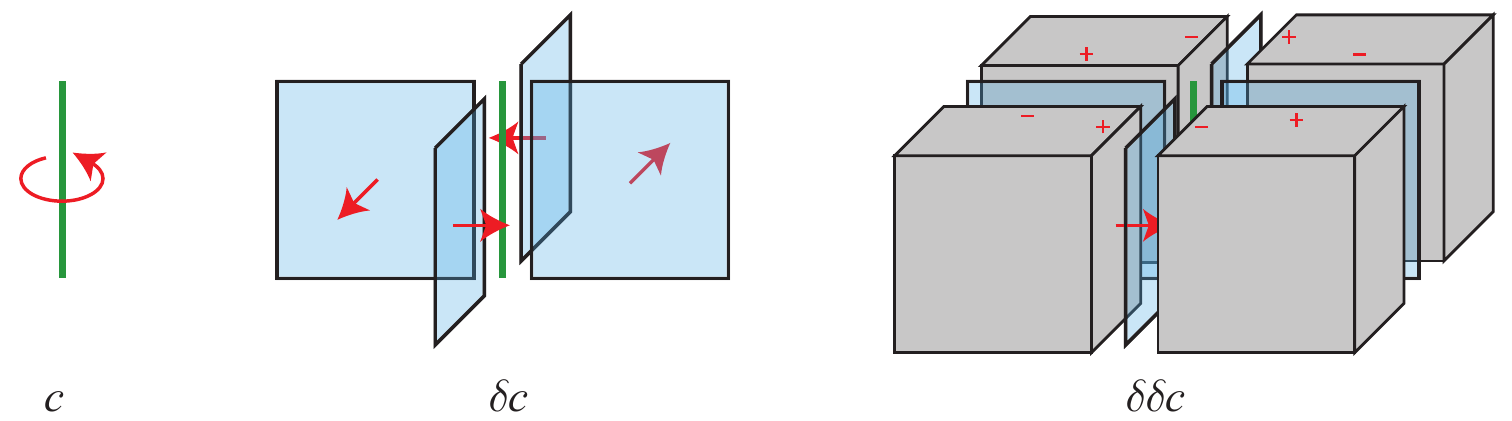}
\caption{The action of twice the coboundary operator $\delta$ on a 1-cell has a zero net result on its surrounding 3-cells, because they all have both a positive and a negative contribution from its neighboring 2-cells (reproduced from \cite{bochevhyman2006}).}
\label{fig:coboundarycoboundary}
\end{figure}

Let $\{\tau_{(k),i}\}$ form a basis of $C_k(D)$, then there is a dual basis $\{\tau^{(k),i}\}$ of $C^k(D)$, such that $\tau^{(k),i}(\tau_{(k),i})=\delta^i_j$ and all $k$-cochains can be represented as linear combinations of the basis elements,
\begin{equation}
\kcochain{c}{k}=\sum_ic_i\tau^{(k),i}.
\end{equation}
The cochains are the discrete analogue of differential forms. With this duality relation between chains and cochains, we can define the formal adjoint of the boundary operator which constitutes an exact sequence on the spaces of $k$-cochains in the cell complex. This formal adjoint is called the \emph{coboundary operator}, $\delta:\kcochainspacedomain{k}{D}\rightarrow\kcochainspacedomain{k+1}{D}$, and is defined as
\begin{equation}
\duality{\delta\kcochain{c}{k}}{\kchain{c}{k+1}} := \duality{\kcochain{c}{k}}{\partial\kchain{c}{k+1}}, \quad\forall\kcochain{c}{k}\in\kcochainspacedomain{k}{D} \text{ and  } \,\forall\kchain{c}{k+1}\in\kchainspacedomain{k+1}{D} \;. \label{algTop_codifferential_dual}
\end{equation}
Also the coboundary operator satisfies $\delta\delta\kcochain{c}{k}=\kcochain{0}{k+2}$ for all $\kcochain{c}{k}\in C^k(D)$, see \figref{fig:coboundarycoboundary}, and gives rise to an exact sequence, called the \emph{cochain complex} $(C^k(D),\delta)$,
\begin{equation}
\label{cochaincomplex}
\begin{CD}
\cdots @>\delta>> C^{k-1}(D)@>\delta>> C^k(D)@>\delta>>C^{k+1}(D)@>\delta>>\cdots\;.
\end{CD}
\end{equation}
The coboundary operator is the discrete analogue of the exterior derivative. Also the coboundary operator has a matrix representation. As a result of the duality pairing in \eqref{algTop_codifferential_dual}, the matrix representation of the coboundary operator is the transpose of the incidence matrix of the boundary operator, $\incidencederivative{k}{k-1}:=\left(\incidenceboundary{k-1}{k}\right)^T$. And again, $\incidencederivative{k+1}{k}\incidencederivative{k}{k-1}=0$. Note that expression \eqref{algTop_codifferential_dual} is nothing but a discrete generalized Stokes' theorem. The matrices representing the coboundary operator only depend on the mesh topology. These matrices will explicitly appear in the final matrix system, \eqref{matrixsystem}.

\subsection{Mimetic Operators} \label{mimeticoperators}
The discretization of the flow variables involves a projection operator, $\pi_h$, from the complete space $\Lambda^k(\Omega)$ to a subspace $\Lambda^k_h(\Omega;C_k)\subset\Lambda^k(\Omega)$. In this subspace differential forms are expressed in terms of $k$-cochains defined on $k$-chains, and corresponding $k$-form interpolation functions (often called basis-functions). Usually, the subspace is a polynomial space. The projection operation actually consists of two steps, a reduction operator, $\reduction$, that integrates the $k$-forms on $k$-chains to get $k$-cochains, and a reconstruction operator, $\reconstruction$, to reconstruct $k$-forms from $k$-cochains using the appropriate basis-functions. These mimetic operators were already introduced before in \cite{bochevhyman2006,HymanScovel1988}. A composition of the two operators gives the projection operator $\projection=\reconstruction\circ\reduction$ as is illustrated below.
\begin{diagram}
\Lambda^k(\Omega)&\rTo^{\quad\projection\quad}& \Lambda^k_h(\Omega;C_k)\\
\dTo^{\reduction} & \ruTo_{\reconstruction} & \\
C^k(D) & &
\end{diagram}
We already saw the similarities between differential geometry and algebraic topology. We now impose constraints on the maps $\reduction$ and $\reconstruction$ to ensure that these structures are preserved.
By imposing structure-preserving constraints on these operations, these three operators together set up the mimetic framework. An extensive discussion on mimetic operators can be found in \cite{kreeftpalhagerritsma2011}. Here only the most important properties are listed.
\begin{definition}[\textbf{Reduction}]
\label{def:reduction}
The {\it reduction operator} $\reduction:\kformspacedomain{k}{\Omega}\rightarrow \kcochainspacedomain{k}{D}$ maps differential forms to cochains. This map is defined by integration as
\begin{equation}
\duality{\reduction \kdifform{a}{k}}{\tau_{(k)}}:=\int_{\tau_{(k)}}\kdifform{a}{k},\quad \forall\tau_{(k)}\in C_k(D).
\label{reduction}
\end{equation}
Then for all $\kchain{c}{k}\in C_k(D)$, the reduction of the $k$-form, $\kdifform{a}{k}\in\Lambda^{k}(\Omega)$, to the $k$-cochain, $\kcochain{a}{k}\in C^k(D)$, is given by
\begin{equation}
\kcochain{a}{k}(\kchain{c}{k}):=\duality{\reduction\kdifform{a}{k}}{\kchain{c}{k}}\stackrel{\eqref{kchain}}{=}\sum_ic^i\duality{\reduction \kdifform{a}{k}}{\tau_{(k),i}}\stackrel{\eqref{reduction}}{=}\sum_ic^i\int_{\tau_{(k),i}}\kdifform{a}{k}=\int_{\kchain{c}{k}}\kdifform{a}{k}.
\end{equation}
\end{definition}
\noindent
The reduction map $\reduction$ provides the {\em integral quantities} that were mentioned in the Introduction.
It is the integration of a $k$-form over all $k$-cells in a $k$-chain that results in a $k$-cochain. A special case of reduction is integration of an $n$-form $a\in\Lambda^n(\Omega)$ over $\Omega$, then
\[
\int_\Omega\kdifform{a}{n}:=\duality{\reduction\kdifform{a}{n}}{\boldsymbol\sigma_{(n)}}\;,
\]
where the chain $\boldsymbol\sigma_{(n)}=\sum_i\tau_{(n),i}$ (so all $c^i=+1$) covers the entire computational domain $\Omega$. The reduction map has a commuting property with respect to continuous and discrete differentiation,
\begin{equation}
\reduction\ederiv=\dederiv\reduction\quad\mathrm{on}\ \Lambda^k(\Omega).
\label{cdp1}
\end{equation}
This commutation can be illustrated as
\[\begin{CD}
\Lambda^k @>\ederiv>> \Lambda^{k+1}\\
@VV\reduction V @VV\reduction V  \\
C^k @>\delta>> C^{k+1}
\end{CD}\]
This property follows from the generalized Stokes' theorem \eqref{stokestheorem} and the duality pairing of \eqref{algTop_codifferential_dual},
\begin{equation}
\langle\reduction\ederiv \kdifform{a}{k},\kchain{c}{k}\rangle\stackrel{\eqref{reduction}}{=}\int_{\kchain{c}{k}}\ederiv \kdifform{a}{k}\stackrel{\eqref{stokestheorem}}{=}\int_{\partial\kchain{c}{k}}\kdifform{a}{k}\stackrel{\eqref{reduction}}{=}\langle\reduction \kdifform{a}{k},\partial\kchain{c}{k}\rangle\stackrel{\eqref{algTop_codifferential_dual}}{=}\langle\dederiv\reduction \kdifform{a}{k},\kchain{c}{k}\rangle.
\end{equation}
The operator acting in the opposite direction to the reduction operator is the reconstruction operator, $\reconstruction$. The \emph{reconstruction operator} $\reconstruction:\kcochainspacedomain{k}{D}\rightarrow\kformspace{k}_h(\Omega;C_k)$ maps $k$-cochains onto finite dimensional $k$-forms. The reconstructed differential forms belong to the space $\Lambda^k_h(\Omega;C_k)$, which is a proper subset of the complete $k$-form space $\Lambda^k(\Omega)$. While the reduction step is clearly defined in \defref{def:reduction}, in the choice of interpolation forms there exists some freedom.
\begin{definition}[\textbf{Reconstruction}] \label{reconstructionoperator}
Although the choice of a reconstruction method allows for some freedom, $\reconstruction$ must satisfy the following properties:
\begin{itemize}
\item Reconstruction $\reconstruction$ must be the right inverse of $\reduction$, so it returns identity ({\it consistency property}),
\begin{equation}
\reduction\reconstruction=Id\quad\mathrm{on}\ C^k(D).
\label{consistency}
\end{equation}
\item Like $\reduction$, also the reconstruction operator $\reconstruction$ has to possess a commuting property with respect to differentiation. A properly chosen reconstruction operator $\reconstruction$ must satisfy a commuting property with respect to the exterior derivative and coboundary operator,
\begin{equation}
\ederiv\reconstruction=\mathcal{I}\dederiv\quad\mathrm{on}\ C^k(D).
\label{cdp2}
\end{equation}
This commutation can be illustrated as
\[\begin{CD}
\kformspaceh{k} @>\ederiv>> \kformspaceh{k+1}\\
@AA\reconstruction A @AA\reconstruction A  \\
\kcochainspace{k} @>\delta>> \kcochainspace{k+1}
\end{CD}\]
\end{itemize}
\end{definition}
Moreover, we want it to be an approximate left inverse of $\reduction$, so the result is close to identity ({\it approximation property})
\begin{equation}
\reconstruction\reduction=Id+\mathcal{O}\left(h^p\right)\quad \mathrm{in}\ \Lambda^k(\Omega).
\label{approximation}
\end{equation}
where $\mathcal{O}(h^p)$ indicates a truncation error in terms of a measure of the mesh size, $h$, and a polynomial order $p$.
\begin{definition}[\textbf{Projection}]\label{th:projection}
The composition $\reconstruction\circ\reduction$ will denote the projection operator, $\projection\define\reconstruction\reduction:\kformspace{k}(\Omega)\rightarrow\kformspace{k}_h(\Omega;C_k)$, allowing for an approximate continuous representation of a $k$-form $a^{(k)}\in\Lambda^k(\Omega)$,
\begin{equation}
\kdifformh{a}{k}=\projection\kdifform{a}{k}=\reconstruction\reduction\kdifform{a}{k}, \quad \pi_h\kdifform{a}{k}\in\kformspace{k}_h(\Omega;C_k)\subset\kformspace{k}(\Omega).
\label{projection}
\end{equation}
where $\reconstruction\reduction\kdifform{a}{k}$ is expressed as a combination of $k$-cochains and interpolating $k$-forms.
\end{definition}
A proof that $\pi_h$ is indeed a projection operator is given in \cite{kreeftpalhagerritsma2011}. Since $\projection\kdifform{a}{k}=\reconstruction\reduction\kdifform{a}{k}$ is a linear combination of $k$-cochains and interpolation $k$-forms, the expansion coefficients in the spectral element formulation to be discussed in Section~\ref{sec:mimeticsem} are the cochains which in turn are the integral quantities mentioned in the Introduction.
\begin{lemma}[\textbf{Commutation property}]
\label{Lem:projectionextder}
There exists a commuting property for the projection and the exterior derivative, such that
\begin{equation}
\label{projectionextder}
\ederiv\projection=\projection\ederiv\quad\mathrm{on}\ \Lambda^k(\Omega).
\end{equation}
This can be illustrated as
\[
\begin{CD}
\Lambda^k @>\ederiv>> \Lambda^{k+1}\\
@VV\projection V @VV\projection V\\
\Lambda_h^k @>\ederiv>> \Lambda_h^{k+1}.
\end{CD}
\]
\begin{proof}
This is a direct consequence of the definitions of the reduction \eqref{cdp1}, reconstruction \eqref{cdp2} and projection operators \eqref{projection},
\[
\ederiv\projection\kdifform{a}{k}\stackrel{\eqref{projection}}{=}\ederiv\reconstruction\reduction\kdifform{a}{k}\stackrel{\eqref{cdp2}}{=}\reconstruction\delta\reduction\kdifform{a}{k}\stackrel{\eqref{cdp1}}{=}\reconstruction\reduction\ederiv\kdifform{a}{k}\stackrel{\eqref{projection}}{=}\projection\ederiv\kdifform{a}{k},\quad\forall \kdifform{a}{k}\in\Lambda^k(\Omega).
\]
\end{proof}
\end{lemma}
Note that it is the intermediate step $\reconstruction\delta\reduction\kdifform{a}{k}$ that is used in practice for the discretization, see Examples \ref{ex:curl} and \ref{ex:div}, \secref{pointwisedivergencefree}. Since we have a matrix representation of the coboundary operator in terms of incidence matrices, we expect the incidence matrices to appear explicitly in the spectral element formulation, see \eqref{matrixsystem}.
\lemmaref{Lem:projectionextder} is the most important result in this paper. As a direct consequence we obtain the pointwise divergence-free solution, as illustrated in the following example.
\begin{example}
Consider the relation $\ud\kdifform{u}{n-1}=\kdifform{g}{n}$. In vector notation the $\ud$ represents the $\mathrm{div}$ operator. Now let $\ud\kdifformh{u}{n-1}=\kdifformh{g}{n}$ be the discretization of our continuous problem. Then by using \eqref{projectionextder} we get
\[
\ud \kdifformh{u}{n-1}-\kdifformh{g}{n}=\ud\pi_h\kdifform{u}{n-1}-\pi_h\kdifform{g}{n}=\pi_h(\ud \kdifform{u}{n-1}-\kdifform{g}{n})=0.
\]
It follows that our discretization is exact. In case $\kdifform{g}{n}=0$, we have a pointwise divergence-free solution of $\kdifformh{u}{n-1}$.
\end{example}
As a direct consequence of \lemmaref{Lem:projectionextder} we satisfy the LBB stability criteria, see \cite{brezzifortin,giraultraviart,kreefterrorestimate}. The projection does \emph{not} commute with codifferential operator. This is the main reason why we rewrite the codifferentials into exterior derivatives and boundary integrals, by means of integration by parts using \eqref{integrationbyparts}.

We do not restrict ourselves to affine mappings only, as is required for many other compatible finite elements, like N\'ed\'elec and Raviart-Thomas elements and their generalizations \cite{arnoldfalkwinther2006,nedelec1980,raviartthomas1977}, but also allow non-affine maps such as transfinite mappings \cite{gordonhall1973} or isogeometric transformations. This allows for better approximations in complex domains with curved boundaries, without the need for excessive refinement. This is possible since the projection operator $\projection$ commutes with the pullback $\pullback$,
\begin{equation}
\label{pullbackprojection}
\pullback\projection=\projection\pullback\quad\mathrm{on}\ \Lambda^k(\Omega).
\end{equation}
This commutation can be illustrated as
\[
\begin{CD}
\Lambda^k(\Omega) @>\pullback>> \Lambda^k(\Omega_{\rm ref})\\
@VV\projection V @VV\projection V\\
\Lambda^k_h(\Omega,C_k) @>\pullback>> \Lambda^k_h(\Omega_{\rm ref},C_k)
\end{CD}
\]
An extensive proof is given in \cite{kreeftpalhagerritsma2011}.

\subsection{Mimetic spectral element basis-functions}\label{sec:mimeticsem}
Now that a mimetic framework is formulated using differential geometry, algebraic topology and the relations between those - the mimetic operators - we derive reconstruction functions, $\mathcal{I}$, that satisfy the properties of the mimetic operators. In combination with the reduction operator, $\mathcal{R}$, it defines the mimetic projection operators, $\pi_h$. The finite dimensional $k$-forms used in this paper are polynomials, based on the idea of spectral element methods, \cite{canuto1}. Spectral element methods have many desirable features such as arbitrary polynomial representation, favourable conditioning, element wise local support,  and optimal stability and approximation properties. However, the definition of the reconstruction operator requires a new set of spectral element interpolation functions. The \emph{mimetic spectral elements} were derived independently by \cite{gerritsma2011,robidoux2008}, and are more extensively discussed in \cite{kreeftpalhagerritsma2011}. Only the most important properties of the mimetic spectral element method are presented here.

In spectral element methods the computational domain $\Omega$ is decomposed into $M$ non-overlapping, possibly curvilinear quadrilateral or hexahedral, closed sub-domains $Q_m$,
\begin{equation}
\Omega=\bigcup_{m=1}^MQ_m,\quad Q_m\cap Q_l=\p Q_m\cap\p Q_l,\ m\neq l,
\end{equation}
where in each sub-domain a Gauss-Lobatto mesh is constructed, see Figures \ref{fig:meshes} and \ref{fig:LDC} in the next section. The complete mesh is indicated by $\mathcal{Q}:=\sum_{m=1}^MQ_m$.

The collection of Gauss-Lobatto meshes in all elements $Q_m\in\mathcal{Q}$ constitutes the cell complex $D$. For each element $Q_m$ there exists a sub cell complex, $D_m$. Note that $D_m\cap D_l,\ m\neq l$, is not an empty set in case they are neighboring elements, but contains all $k$-cells, $k<n$, of the common boundary, see \defref{cellcomplex}.

Each sub-domain is mapped from the reference element, $Q_{\rm ref}=[-1,1]^n$, using the mapping $\Phi_m:Q_{\rm ref}\rightarrow Q_m$. Then all flow variables defined on $Q_m$ are pulled back onto this reference element using the following pullback operation, $\Phi^\star_m:\Lambda^k_h(Q_m;C_k)\rightarrow\Lambda^k_h(Q_{\rm ref};C_k)$. In three dimensions the reference element is given by $Q_{\rm ref}:=\{(\xi,\eta,\zeta)\;|\;-1\leq\xi,\eta,\zeta\leq1\}$.

The basis-functions that interpolate the cochains on the quadrilateral or hexahedral elements are constructed using tensor products. It is therefore sufficient to derive interpolation functions in one dimension and use tensor products afterwards to construct $n$-dimensional basis functions. A similar approach was taken in \cite{buffa2011b}. Because the projection operator and the pullback operator commute \eqref{pullbackprojection}, the interpolation functions are discussed for the reference element only.

\begin{figure}[t]
      \begin{minipage}[t]{0.49\linewidth}
            \centering\includegraphics[width=0.9\linewidth]{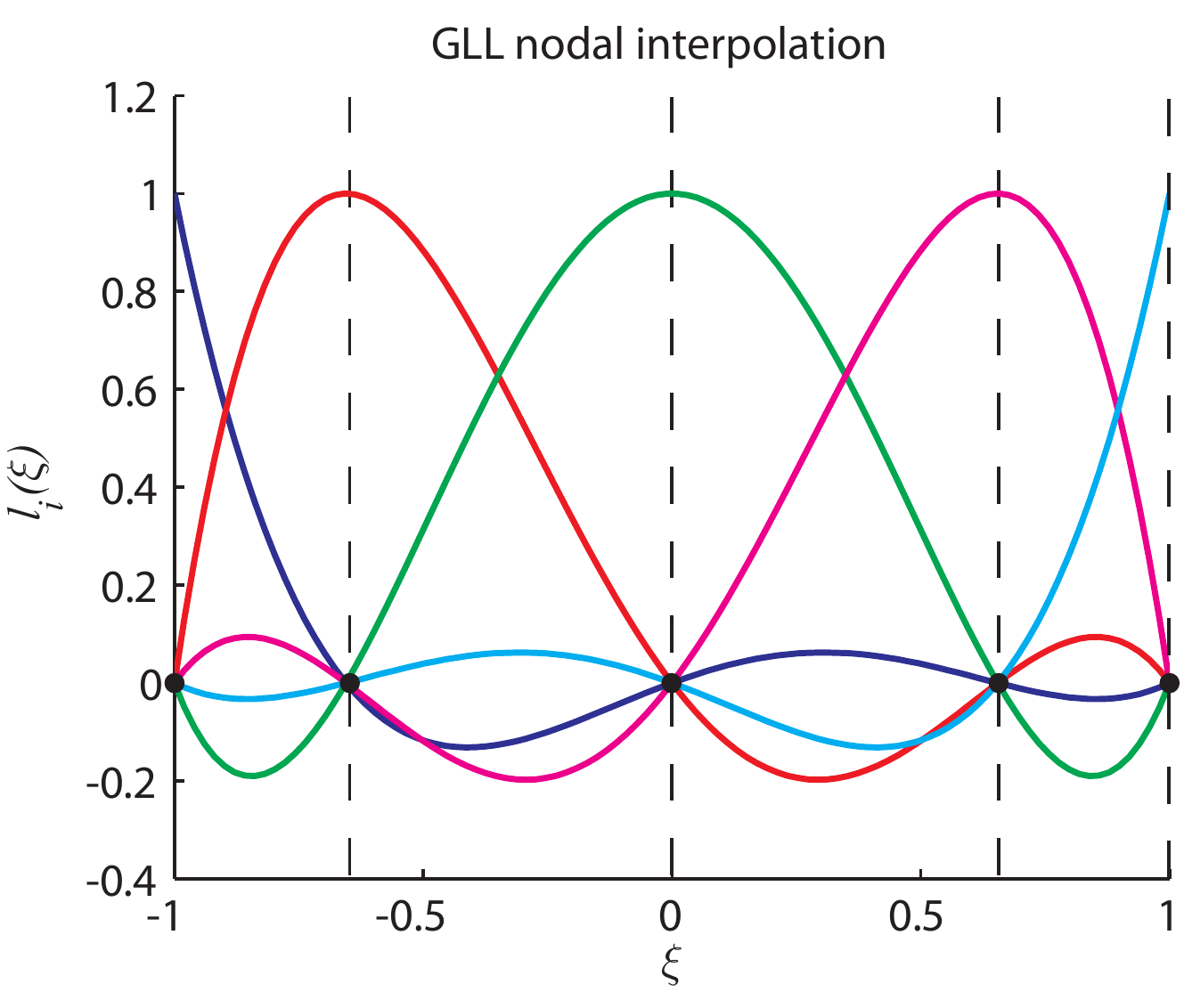}
	\caption{Lagrange polynomials on Gauss-Lobatto-Legendre mesh.}
	\label{fig:lagrange}
    \end{minipage}\hfill
    \begin{minipage}[t]{0.49\linewidth}
            \centering\includegraphics[width=0.9\linewidth]{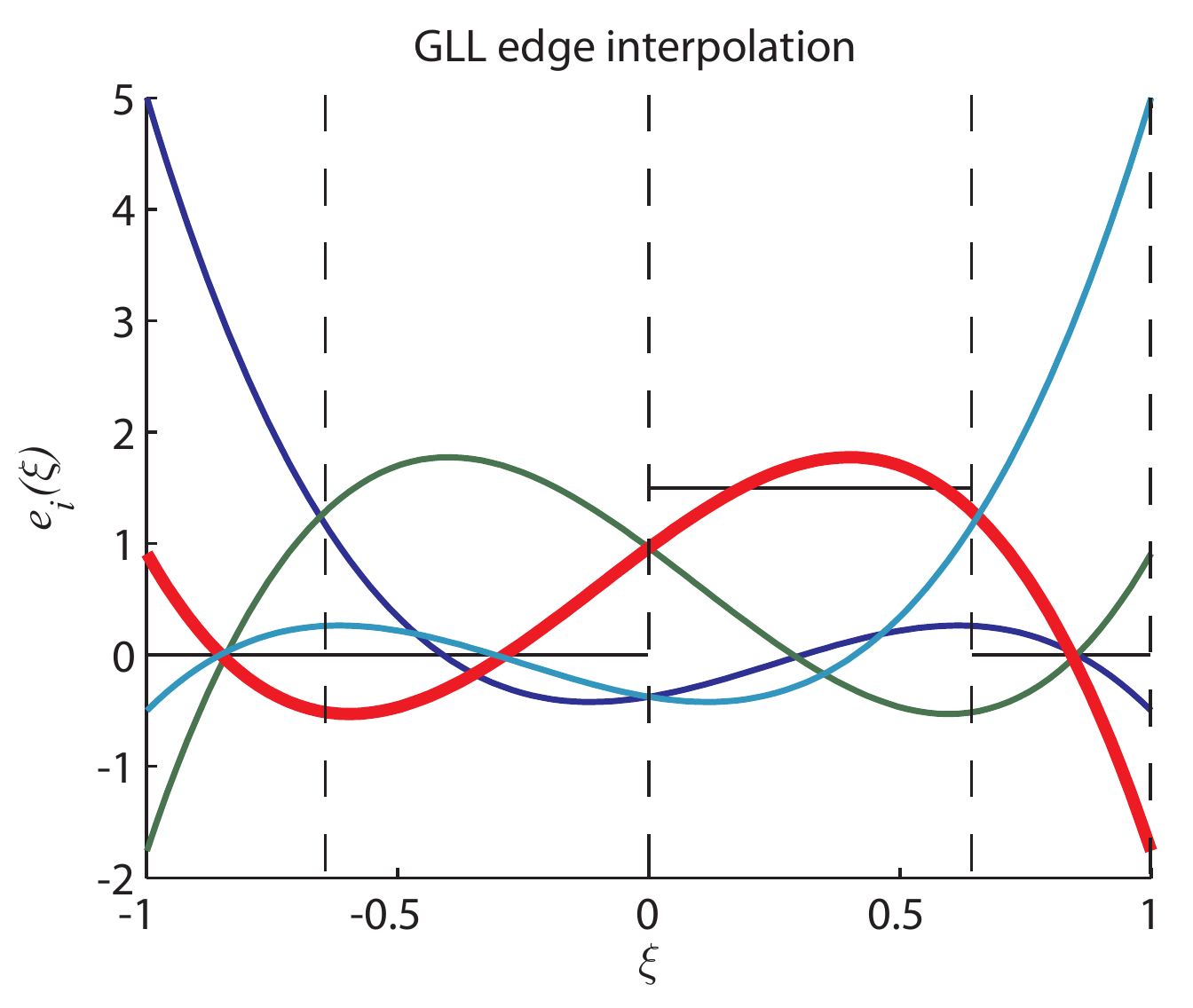}
	\caption{Edge polynomials on Gauss-Lobatto-Legendre mesh.}
	\label{fig:edgepoly}
    \end{minipage}
\end{figure}

Consider a 0-form $\kdifform{a}{0}\in\Lambda^0(Q_{\rm ref})$ on $Q_{\rm ref}:=\xi\in[-1,1]$, on which a cell complex $D$ is defined that consists of $N+1$ nodes, $\xi_i$, where $-1\leq \xi_0<\hdots<\xi_N\leq 1$, and $N$ edges, $\tau_{(1),i}=[\xi_{i-1},\xi_i]$, of which the nodes are the boundaries. Corresponding to this set of nodes (0-chains) there exists a projection using $N^{\rm th}$ order \textit{Lagrange polynomials}, $l_i(\xi)$, to approximate a $0$-form, as
\begin{equation}
\projection\kdifform{a}{0}=\sum_{i=0}^N a_il_i(\xi).
\label{nodalapprox}
\end{equation}
Lagrange polynomials have the property that they interpolate nodal values and are therefore suitable to reconstruct the cochain $\kcochain{a}{0}=\reduction\kdifform{a}{0}$ containing the set $a_i=a(\xi_i)$ for $i=0,\hdots,N$. So Lagrange polynomials can be used to reconstruct a 0-form from a 0-cochain. Lagrange polynomials are in fact 0-forms themselves, $l_i(\xi)\in\Lambda^0_h(Q_{\rm ref};C_0)$.
Lagrange polynomials are constructed such that their value is one in the corresponding point and zero in all other mesh points,
\begin{equation}
 \reduction l_i(\xi)=l_i(\xi_p)=\left\{
\begin{aligned}
&1&{\rm if}\ i=p\\ &0&{\rm if}\ i\neq p
\end{aligned}
\right..
 \label{nodalproperty}
\end{equation}
This satisfies \eqref{consistency}, where in this case $\reconstruction=l_i(\xi)$. Gerritsma \cite{gerritsma2011} and Robidoux \cite{robidoux2008} derived a similar projection for 1-forms, consisting of $1$-cochains and $1$-form polynomials, that is called the \textit{edge polynomial}, $e_i(\xi)\in\Lambda^1_h(Q_{\rm ref};C_1)$.
\begin{lemma}[\textbf{Edge polynomial}]
\label{lemma:edge}
Following Definitions \ref{def:reduction} and \ref{reconstructionoperator}, apply the exterior derivative to $\pi_h\kdifform{a}{0}$, it gives the 1-form $\pi_h\kdifform{b}{1}=\ederiv \pi_h\kdifform{a}{0}=\reconstruction\delta\reduction \kdifform{a}{0}$ given by
\begin{equation}
\projection \kdifform{b}{1}=\sum_{i=1}^Nb_i e_i(\xi),
\end{equation}
with 1-cochain $\kcochain{b}{1}$, where
\begin{align}
\label{ucochain}
b_i&=\langle\reduction \kdifform{b}{1},\tau_{(1),i}\rangle=\int_{\tau_{(1),i}}b(\xi)=\int_{\tau_{(1),i}}\ederiv \kdifform{a}{0}=\int_{\p\tau_{(1),i}}\kdifform{a}{0},\\
&=a(\xi_i)-a(\xi_{i-1})=a_i-a_{i-1},\nonumber
\end{align}
with the edge interpolation polynomial defined by
\begin{align}
\label{edge}
e_i(\xi)&=-\sum_{k=0}^{i-1}\ederiv l_k(\xi)=\sum_{k=i}^{N}\ederiv l_k(\xi)=\tfrac{1}{2}\sum_{k=i}^{N}\ederiv l_k(\xi)-\tfrac{1}{2}\sum_{k=0}^{i-1}\ederiv l_k(\xi).
\end{align}
\begin{proof}
See \cite{gerritsma2011,kreeftpalhagerritsma2011,robidoux2008}.
\end{proof}
\end{lemma}
The value corresponding to line segment (1-cell) $\tau_{(1),i}$ is given by $b_i=a_i-a_{i-1}$ and so $\kcochain{b}{1}=\delta\kcochain{a}{0}$ is the discrete derivative operator in 1D. This operation is purely topological, no metric is involved. It satisfies \eqref{cdp2}, since $\ederiv\reconstruction\kcochain{a}{0}=\reconstruction\delta\kcochain{a}{0}$. Note that we have $\ederiv e_i(\xi)=\sum\ederiv\!\circ\!\ederiv  l_i(\xi)=0$.
The 1-form edge polynomial can also be written as below, separating the edge function into its polynomial and its basis,
\begin{displaymath}
e_i(\xi)=\ve_i(\xi)\ederiv \xi,\quad\mathrm{with}\quad\ve_i(\xi)=-\sum_{k=0}^{i-1}\frac{\ederiv l_k}{\ederiv \xi}.
\end{displaymath}
Similar to \eqref{nodalproperty}, the edge functions are constructed such that when integrating $e_i(\xi)$ over a line segment it gives one for the corresponding element and zero for any other line segment, so
\begin{equation}
 \reduction e_i(\xi)=\int_{\xi_{p-1}}^{\xi_p}e_i(\xi)=\left\{
\begin{aligned}
&1&{\rm if}\ i=p\\ &0&{\rm if}\ i\neq p
\end{aligned}
\right..
 \label{intedge}
\end{equation}
This also satisfies \eqref{consistency}, where in this case $\reconstruction=e_i(\xi)$.
The fourth-order Lagrange and third-order edge polynomials, corresponding to a Gauss-Lobatto mesh with $N=4$, are shown in Figures \ref{fig:lagrange} and \ref{fig:edgepoly}.

Now that we have developed interpolation functions in one dimension, we can extend this to the multidimensional framework by means of tensor products. This allows for the interpolation of integral quantities defined on $k$-dimensional cubes. Consider a reference element in $\mathbb{R}^3$, $Q_{\rm ref}=[-1,1]^3$. Then the interpolation functions for points, lines, surfaces and volumes are given by,
\begin{displaymath}
\begin{aligned}
&\mathrm{point}:&&P^{(0)}_{i,j,k}(\xi,\eta,\zeta)=l_i(\xi)\otimes l_j(\eta)\otimes l_k(\zeta),\\
&\mathrm{line}:&&L^{(1)}_{i,j,k}(\xi,\eta,\zeta)=\{e_i(\xi)\otimes l_j(\eta)\otimes l_k(\zeta),\ l_i(\xi)\otimes e_j(\eta)\otimes l_k(\zeta),\ l_i(\xi)\otimes l_j(\eta)\otimes e_k(\zeta)\},\\
&\mathrm{surface}:&&S^{(2)}_{i,j,k}(\xi,\eta,\zeta)=\{l_i(\xi)\otimes e_j(\eta)\otimes e_k(\zeta),\ e_i(\xi)\otimes l_j(\eta)\otimes e_k(\zeta),\ e_i(\xi)\otimes e_j(\eta)\otimes l_k(\zeta)\},\\
&\mathrm{volume}:&&V^{(3)}_{i,j,k}(\xi,\eta,\zeta)=e_i(\xi)\otimes e_j(\eta)\otimes e_k(\zeta).
\end{aligned}
\end{displaymath}
Note that $V^{(3)}_{i,j,k}$ is indeed a 3-form, since $e_i(\xi)\otimes e_j(\eta)\otimes e_k(\zeta)=\ve_i(\xi)\ve_j(\eta)\ve_k(\zeta)\,\ud\xi\wedge\ud\eta\wedge\ud\zeta$. So the approximation spaces are spanned by combinations of Lagrange and edge basis functions,
\begin{align*}
\Lambda^0_h(\mathcal{Q};C_0)&:=\mathrm{span}\left\{P^{(0)}_{i,j,k}\right\}_{i=0,j=0,k=0}^{N,N,N},\\
\Lambda^1_h(\mathcal{Q};C_1)&:=\mathrm{span}\left\{\big(L^{(1)}_{i,j,k}\big)_1\right\}_{i=1,j=0,k=0}^{N,N,N}\times \mathrm{span}\left\{\big(L^{(1)}_{i,j,k}\big)_2\right\}_{i=0,j=1,k=0}^{N,N,N}\times \mathrm{span}\left\{\big(L^{(1)}_{i,j,k}\big)_3\right\}_{i=0,j=0,k=1}^{N,N,N},\\
\Lambda^2_h(\mathcal{Q};C_2)&:=\mathrm{span}\left\{\big(S^{(2)}_{i,j,k}\big)_1\right\}_{i=0,j=1,k=1}^{N,N,N}\times \mathrm{span}\left\{\big(S^{(2)}_{i,j,k}\big)_2\right\}_{i=1,j=0,k=1}^{N,N,N}\times \mathrm{span}\left\{\big(S^{(2)}_{i,j,k}\big)_3\right\}_{i=1,j=1,k=0}^{N,N,N},\\
\Lambda^3_h(\mathcal{Q};C_3)&:=\mathrm{span}\left\{V^{(3)}_{i,j,k}\right\}_{i=1,j=1,k=1}^{N,N,N}.
\end{align*}

Lagrange interpolation by itself does not guarantee a convergent approximation \cite{erdos1980}, but it requires a suitably chosen set of points, $-1\leq\xi_0<\xi_1<\hdots<\xi_N\leq1$. Here, the Gauss-Lobatto distribution is proposed, because of its superior convergence behaviour \cite{canuto1}. The convergence rates of Lagrange and edge interpolants were obtained in \cite{kreeftpalhagerritsma2011} and are given by,
\begin{align}
\Vert \kdifform{a}{0}-\pi_h\kdifform{a}{0}\Vert_{H\Lambda^0}&\leq C\frac{h^{l-1}}{p^{m-1}}|\kdifform{a}{0}|_{H^m\Lambda^0},\label{hpconvergenceestimate}\\
\Vert \kdifform{b}{1}-\pi_h\kdifform{b}{1}\Vert_{L^2\Lambda^1}&\leq C\frac{h^{l-1}}{p^{m-1}}|\kdifform{b}{1}|_{H^{m-1}\Lambda^1},\label{edgeinterpolationerror}
\end{align}
with $l=\mathrm{min}(p+1,m)$. For the variables vorticity, velocity and pressure in the VVP formulation of the Stokes problem, the $h$-convergence rates of the interpolation errors become,
\begin{gather}
\label{vvpinterpolationerror}
\norm{\omega-\pi_h\omega}_{L^2\Lambda^{n-2}}=\mathcal{O}(h^{N+s}),\quad\norm{\omega-\pi_h\omega}_{H\Lambda^{n-2}}=\mathcal{O}(h^N),\nonumber\\
\norm{u-\pi_h u}_{L^2\Lambda^{n-1}}=\mathcal{O}(h^N),\quad\norm{p-\pi_h p}_{L^2\Lambda^{n}}=\mathcal{O}(h^N),
\end{gather}
where $s=1$ for $n=2$ and $s=0$ for $n>2$, and with $N$ defined as in \secref{sec:mimeticsem}. Because of \eqref{stokeseq3} and \eqref{projectionextder}, we have $\norm{u-\pi_h u}_{H\Lambda^{n-1}}=\norm{u-\pi_h u}_{L^2\Lambda^{n-1}}$.

\subsection{Pointwise divergence-free discretization}\label{pointwisedivergencefree}
One of the most interesting properties of the mimetic method presented in this paper, is that within our weak formulation, the divergence-free constraint is satisfied pointwise. This result follows from the three commuting properties with the exterior derivative, \eqref{cdp1}, \eqref{cdp2} and \eqref{projectionextder}, as was shown in \lemmaref{Lem:projectionextder}. The corresponding commuting diagrams are repeated in the diagram below for the two dimensional case.
\begin{diagram}
\Lambda^{0}_h(\mathcal{Q};C_0) & \rTo^{\quad\ederiv\quad}_{\mathrm{curl}} & \Lambda^{1}_h(\mathcal{Q};C_1) & \rTo^{\quad\ederiv\quad}_{\mathrm{div}} & \Lambda^{2}_h(\mathcal{Q};C_2)\\
\dTo^{\reduction} \uTo_{\reconstruction} & & \dTo^{\reduction} \uTo_{\reconstruction} &  & \dTo^{\reduction} \uTo_{\reconstruction}\\
C^{0}(D) & \rTo^{\delta} & C^{1}(D) & \rTo^{\delta} & C^{2}(D)
\end{diagram}
Note that by curl we refer to the two-dimensional variant, applied to a scalar, i.e. $\mathrm{curl}\,\omega=(\p\omega/\p y,-\p\omega/\p x)^T$, see also Example~\ref{ex:2dstokes}, and is also called the normal gradient operator, $\mathrm{grad}^\perp$, see \cite{palha2010}.

In the following two examples we demonstrate the action of the exterior derivative on vorticity, $\kdifformh{\omega}{0}\in\Lambda_h^{0}(Q_{\rm ref};C_0)$, and on the velocity flux, $\kdifformh{u}{1}\in\Lambda^{1}_h(Q_{\rm ref};C_1)$. Two dimensional reconstruction is based on tensor product construction of the one dimensional reconstruction function introduced above.
\begin{example}[\textbf{Curl operator}]\label{ex:curl}
Consider a flux $\kdifformh{z}{1}\in\Lambda^1_h(Q_{\rm ref};C_1)$ with $C_1$ outer-oriented, and where $\kdifformh{z}{1}=\ud\kdifformh{\omega}{0}$. Then $\kdifformh{\omega}{0}$ is expanded in the reference coordinates $(\xi,\eta)$ as
\begin{equation}
\kdifformh{\omega}{0}=\sum_{i=0}^N\sum_{j=0}^N\omega_{i,j}l_i(\xi)l_j(\eta).
\end{equation}
Apply the exterior derivative in the same way as in \lemmaref{lemma:edge}, it gives
\begin{subequations}
\begin{align}
\kdifformh{z}{1}=\ud\kdifformh{\omega}{0}&=\sum_{i=1}^N\sum_{j=0}^N (\omega_{i,j}-\omega_{i-1,j})e_i(\xi)l_j(\eta)+\sum_{i=0}^N\sum_{j=1}^N(\omega_{i,j}-\omega_{i,j-1}) l_i(\xi)e_j(\eta),\nonumber\\
&=-\sum_{i=1}^N\sum_{j=0}^N z^\eta_{i,j}e_i(\xi)l_j(\eta)+\sum_{i=0}^N\sum_{j=1}^Nz^\xi_{i,j} l_i(\xi)e_j(\eta),
\end{align}
\end{subequations}
where $z^\xi_{i,j}=\omega_{i,j}-\omega_{i,j-1}$, and $z^\eta_{i,j}=\omega_{i-1,j}-\omega_{i,j}$ can be compactly written as $\mathbf{z}^{(1)}=\delta\boldsymbol\omega^{(0)}$, with $\boldsymbol\omega^{(0)}\in C^0(D)$ and $\mathbf{z}^{(1)}\in C^1(D)$, or in matrix notation as $\mathbf{z}=\mathsf{E}^{(1,0)}\boldsymbol\omega$. This relation is exact, coordinate free and invariant under transformations.
\end{example}

\begin{example}[\textbf{Divergence operator}]\label{ex:div}
Let $\kdifformh{u}{1}\in\Lambda^1_h(Q_{\rm ref};C_1)$ be the velocity flux defined as
\begin{equation}
\kdifformh{u}{1}=-\sum_{i=1}^N\sum_{j=0}^Nv_{i,j}e_i(\xi)l(\eta)+\sum_{i=0}^N\sum_{j=1}^Nu_{i,j}l_i(\xi)e_j(\eta).
\end{equation}
Compare this to the velocity flux in Example~\ref{ex:2dstokes}, p.\pageref{ex:2dstokes}. Then the change of mass, $\kdifformh{m}{2}\in\Lambda^2_h(Q_{\rm ref};C_2)$, is equal to the exterior derivative of $\kdifformh{u}{1}$,
\begin{align}
\kdifformh{m}{2}=\ud \kdifformh{u}{1}&=\sum_{i=1}^N\sum_{j=1}^N(u_{i,j}-u_{i-1,j}+v_{i,j}-v_{i,j-1})e_i(\xi)e_j(\eta).\nonumber\\
&=\sum_{i=1}^N\sum_{j=1}^Nm_{i,j}e_i(\xi)e_j(\eta),
\end{align}
where $m_{i,j}=u_{i,j}-u_{i-1,j}+v_{i,j}-v_{i,j-1}$ can be compactly written as $\mathbf{m}^{(2)}=\delta\mathbf{u}^{(1)}$, with $\mathbf{u}^{(1)}\in C^1(D)$ and $\mathbf{m}^{(2)}\in C^2(D)$, or in matrix notation as $\mathbf{m}=\mathsf{E}^{(2,1)}\mathbf{u}$. Note that if the mass production is zero, as in our model problem \eqref{mass}, the incompressibility constraint is already satisfied at discrete/cochain level. Interpolation then results in a pointwise divergence-free solution.
\end{example}

\section{Mixed formulation, boundary conditions and implementation}\label{sec:mixedformulationstokes}
We know how to discretize exactly the metric-free exterior derivative $\ud$ (see \lemmaref{Lem:projectionextder}, \secref{mimeticoperators}, and the examples above), but it is less obvious how to treat the codifferential operator $\ud^*$. Fortunately, the two are directly related using $L^2$-inner products as seen in \eqref{integrationbyparts}. Therefore the derivation of the mixed formulation of the Stokes problem consists of two steps: 1). Multiply equations \eqref{stokeseq1}-\eqref{stokeseq3} by the test functions $\kdifform{\tau}{n-2},\kdifform{v}{n-1},\kdifform{q}{n}$ using $L^2$-inner products. 2). Use integration by parts, as in \eqref{integrationbyparts}, to express the codifferentials in terms of the exterior derivatives and boundary integrals. The resulting mixed formulation of the Stokes problem becomes:\\
 \\
Find $(\kdifform{\omega}{n-2},\kdifform{u}{n-1},\kdifform{p}{n})\in\{H\Lambda^{n-2}\times H\Lambda^{n-1}\times L^2\Lambda^n\}$, given $\kdifform{f}{n}\in L^2\Lambda^{n-1}$, for all $(\kdifform{\tau}{n-2},\kdifform{v}{n-1},\kdifform{q}{n})\in\{H\Lambda^{n-2}\times H\Lambda^{n-1}\times L^2\Lambda^n\}$, such that
\begin{subequations}
\label{mixedstokes}
\begin{align}
\big(\kdifform{\tau}{n-2},\kdifform{\omega}{n-2}\big)_\Omega-\big(\ud\kdifform{\tau}{n-2},\kdifform{u}{n-1}\big)_\Omega&=-\int_{\p\Omega}\tr\kdifform{\tau}{n-2}\wedge\tr\star \kdifform{u}{n-1}, \label{mixedstokes1}\\
\big(\kdifform{v}{n-1},\ud\kdifform{\omega}{n-2}\big)_\Omega+\big(\ud \kdifform{v}{n-1},\kdifform{p}{n}\big)_\Omega&=\big(\kdifform{v}{n-1},\kdifform{f}{n-1}\big)_\Omega +\int_{\p\Omega}\tr \kdifform{v}{n-1}\wedge\tr\star \kdifform{p}{n}, \label{mixedstokes2}\\
\big(\kdifform{q}{n},\ud \kdifform{u}{n-1})_\Omega&=0. \label{mixedstokes3}
\end{align}
\end{subequations}
This mixed formulation is similar to those in \cite{bernardi2006,dubois2002,giraultraviart}. The mixed formulation is well-posed, see \cite{giraultraviart,kreefterrorestimate}. The discrete problem is almost the same as the continuous problem, that is: find $(\kdifformh{\omega}{n-2},\kdifformh{u}{n-1},\kdifformh{p}{n})\in\{\Lambda_h^{n-2}\times \Lambda_h^{n-1}\times \Lambda_h^n\}$, given $\kdifformh{f}{n}\in \Lambda_h^{n-1}$, for all $(\kdifformh{\tau}{n-2},\kdifformh{v}{n-1},\kdifformh{q}{n})\in\{\Lambda_h^{n-2}\times \Lambda_h^{n-1}\times \Lambda_h^n\}$, such that \eqref{mixedstokes1}-\eqref{mixedstokes3} hold. The discrete problem is also well-posed, because every subcomplex of a Hilbert complex is also a Hilbert complex, so if $(H\Lambda,\ud)$ is a Hilbert complex, so is $(\Lambda_h,\ud)$, and the projection operator from $H\Lambda^k(\Omega)$ to $\Lambda^k_h(\Omega;C_k)$ is bounded, see \cite{kreeftpalhagerritsma2011}. A complete proof is given in \cite{kreefterrorestimate}.

System \eqref{stokeseq} needs to be supplemented with boundary conditions on $\p\Omega$. Their exists four possible types of boundary conditions as follows from the boundary integrals in the mixed formulation, \eqref{mixedstokes}. Subdivide the boundary into several parts, $\p\Omega=\bigcup_i\Gamma_i$, where $\Gamma_i\cap\Gamma_j=\emptyset$ for $i\neq j$. Each part of the boundary can have one of the following four boundary conditions: 1. prescribed \emph{velocity} (such as no-slip), 2. \emph{tangential velocity - pressure}, 3. \emph{tangential vorticity - normal velocity}, and 4. \emph{tangential vorticity - pressure} boundary conditions. An overview is given in Table~\ref{tab:boundaryconditions}.

\begin{table}[htbp]
\begin{tabular}{lccc}
Name  & $\quad\quad$ Exterior Calculus $\quad\quad$ & $\quad\quad$ Vector Calculus $\quad\quad$ & $\quad$ Type $\quad$ \\ \hline\hline
 & & & \\
Normal velocity & $\tr\kdifform{u}{n-1}\ \Rightarrow\ \tr\kdifform{v}{n-1}=0$ & $\vec{u}\cdot\vec{n}\ \Rightarrow\ \vec{v}\cdot\vec{n}=0$ & essential \\
tangential velocity & $\tr\star\kdifform{u}{n-1}$ & $\vec{u}\cdot\vec{t}$ & natural \\
 & & & \\ \hline  & & & \\
Tangential velocity & $\tr\star\kdifform{u}{n-1}$ & $\vec{u}\cdot\vec{t}$ & natural \\
pressure & $\tr\star\kdifform{p}{n}$& $p$ & natural \\
 & & & \\ \hline  & & & \\
Tangential vorticity & $\tr\kdifform{\omega}{n-2}\ \Rightarrow\ \tr\kdifform{\tau}{n-2}=0$ & $\vec{\omega}\times\vec{t}\ \Rightarrow\ \vec{\tau}\times\vec{t}=\vec{0}$ & essential \\
normal velocity & $\tr\kdifform{u}{n-1}\ \Rightarrow\ \tr\kdifform{v}{n-1}=0$ & $\vec{u}\cdot\vec{n}\ \Rightarrow\ \vec{v}\cdot\vec{n}=0$ & essential \\
 & & & \\ \hline  & & & \\
Tangential vorticity & $\tr\kdifform{\omega}{n-2}\ \Rightarrow\ \tr\kdifform{\tau}{n-2}=0$ & $\vec{\omega}\times\vec{t}\ \Rightarrow\ \vec{\tau}\times\vec{t}=\vec{0}$ & essential \\
pressure & $\tr\star\kdifform{p}{2}$ & $p$ & natural \\
 & & & \\
\end{tabular}
\caption{Admissible boundary conditions for Stokes flow in vorticity-velocity-pressure formulation.}
\label{tab:boundaryconditions}
\end{table}

From the implementation point of view we would like to mention that the $L^2$ inner products and boundary integrals are evaluated using Gauss-Lobatto quadrature, which is exact for polynomials up to order $2N-1$, \cite{canuto1}. The resulting system matrix is a saddle point system that is given by,
\begin{equation}
\label{matrixsystem}
\begin{bmatrix}
\mathsf{M}^{(n-2)} & \big(\mathsf{E}^{(n-1,n-2)}\big)^T\mathsf{M}^{(n-1)} & \emptyset \\
\mathsf{M}^{(n-1)}\mathsf{E}^{(n-1,n-2)} & \emptyset & \big(\mathsf{E}^{(n,n-1)}\big)^T\mathsf{M}^{(n)} \\
\emptyset & \mathsf{M}^{(n)}\mathsf{E}^{(n,n-1)} &  \emptyset
\end{bmatrix}
\begin{bmatrix}
\boldsymbol\omega \\ \mathbf{u} \\ \mathbf{p}
\end{bmatrix}
=
\begin{bmatrix}
-\mathsf{B}_1(\mathbf{\star u}) \\ \mathsf{M}^{(n-1)}\mathbf{f}^{(n-1)}+\mathsf{B}_2(\mathbf{\star p}) \\ \emptyset
\end{bmatrix}
\end{equation}
The final system matrix is symmetric and only consists of $L^2$ inner product matrices for $k$-forms, $\mathsf{M}^{(k)}$ (also known as mass matrices), and incidence matrices, $\mathsf{E}^{(k,k-1)}$, that are directly obtained from the mesh topology, see p.\pageref{cochaincomplex}. Coordinate transformations imposed by the pullback operator appear in the $L^2$ inner products as a standard change of basis, see also \cite{bouman2011}. The matrices $\mathsf{B}_1$ and $\mathsf{B}_2$ represent the boundary integrals in \eqref{mixedstokes1} and \eqref{mixedstokes2}, and $(\star\mathbf{u})$ and $(\star\mathbf{p})$ are the tangential velocity and pressure boundary conditions imposed. A discussion on efficient solvers for symmetric indefinite systems that follow from saddle point problems can be found in \cite{benzigolubliesen2005,rehman2011}.

\section{Numerical Results}\label{sec:numericalresults}
Now that all parts of the mixed mimetic method are treated, we can test the performance of the numerical scheme using a set of three test problems. The first one consists of an analytic solution on a unit square, where optimal $h$-convergence and exponential $p$-convergence rates are shown for both Cartesian and curvilinear meshes for all combinations of boundary conditions. The second is a lid-driven cavity flow, where results are compared with a reference solution. Finally, Stokes flow around a cylinder moving with constant velocity in a channel is considered.

\subsection{Manufactured solution}
\begin{figure}[htb]
\centering
\includegraphics[width=0.7\textwidth]{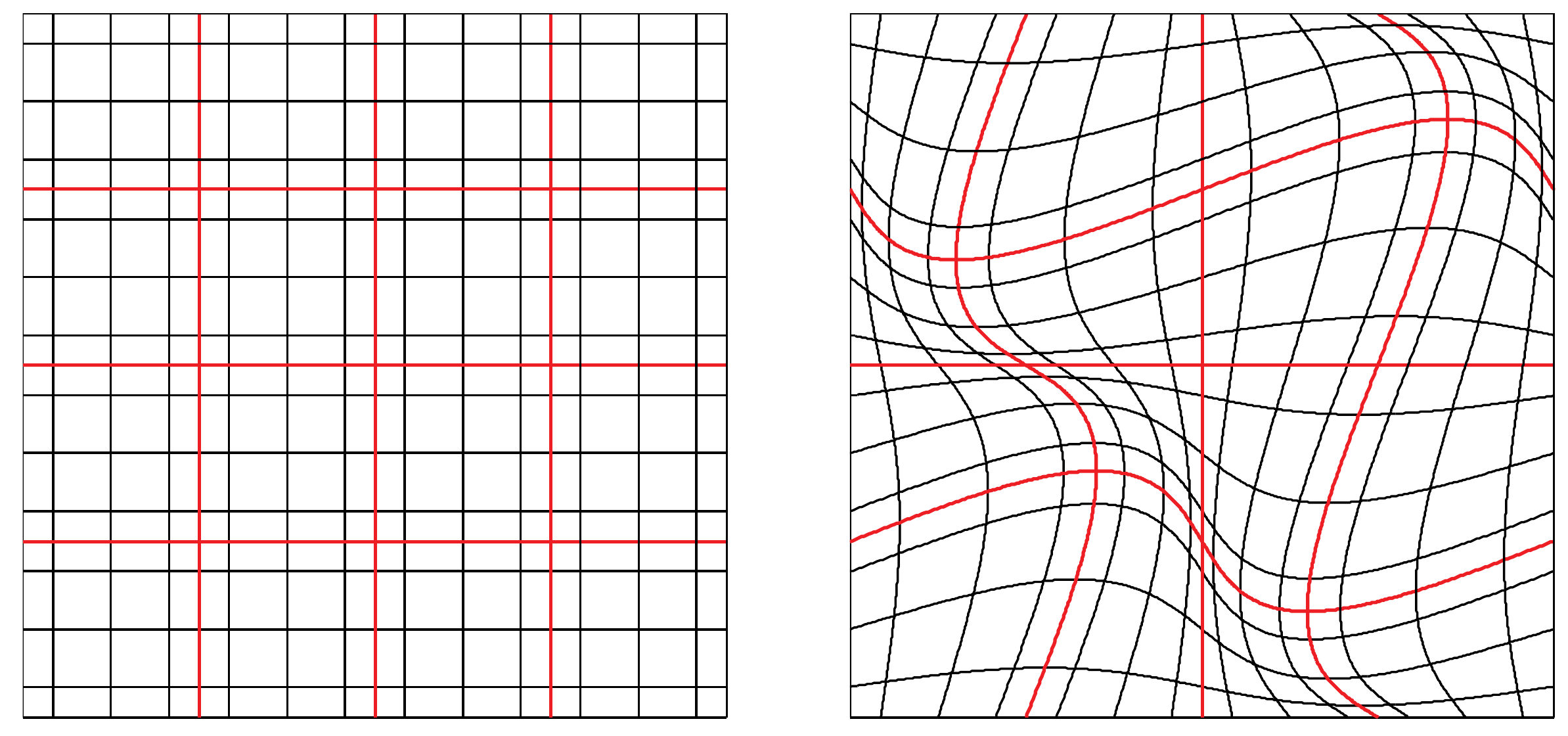}
\caption{Examples of a Cartesian and a curvilinear mesh used in the convergence analysis. The meshes shown consist of $4\times4$ spectral elements, with for each element, $N=4$. The element boundaries are indicated in red.}
\label{fig:meshes}
\end{figure}

The first test case addresses the convergence for $h$- and $p$-refinement of the mixed mimetic spectral element method applied to the Stokes model. The model problem is defined on the unit square $\Omega=[0,1]^2$, with Cartesian coordinates $\mathbf{x}:=(x,y)$, with $\nu=1$ and with the right hand side $\kdifform{f}{1}\in\Lambda^1(\Omega)$ given by
\begin{subequations}
\label{testcase1}
\begin{align}
\kdifform{f}{1}=-&f_y(\mathbf{x})\,\ud x+f_x(\mathbf{x})\,\ud y,\nonumber\\
=-&\left(\pi\sin(\pi x)\cos(\pi y)+8\pi^2\cos(2\pi x)\sin(2\pi y)\right)\ud x\nonumber\\
+&\left(\pi\cos(\pi x)\sin(\pi y)-8\pi^2\sin(2\pi x)\cos(2\pi y)\right)\ud y.
\end{align}
This right hand side results in an exact solution for the vorticity $\kdifform{\omega}{0}\in\Lambda^0(\Omega)$, velocity flux $\kdifform{u}{1}\in\Lambda^1(\Omega)$, and pressure $\kdifform{p}{2}\in\Lambda^2(\Omega)$ components of the Stokes problem, given by
\begin{align}
\kdifform{\omega}{0}&=\omega(\mathbf{x})=-4\pi\sin(2\pi x)\sin(2\pi y),\\
\kdifform{u}{1}&=-v(\mathbf{x})\,\ud x+u(\mathbf{x})\,\ud y\nonumber\\
 &=-\left(\cos(2\pi x)\sin(2\pi y)\right)\ud x+\left(-\sin(2\pi x)\cos(2\pi y)\right)\ud y,\\
\kdifform{p}{2}&=p(\mathbf{x})\,\ud x\!\wedge\!\ud y=\left(\sin(\pi x)\sin(\pi y)\right)\ud x\!\wedge\!\ud y.
\end{align}
\end{subequations}
This testcase was discussed before in \cite{gerritsmaphillips,prootgerritsma2002}. Calculations were performed on both a Cartesian as well as a curvilinear mesh as shown in \figref{fig:meshes}. The map, $(x,y)=\Phi(\xi,\eta)$, used for the curved mesh is given by
\begin{subequations}
\begin{align}
x(\xi,\eta)&=\tfrac{1}{2}+\tfrac{1}{2}\left(\xi+\tfrac{1}{5}\sin(\pi\xi)\sin(\pi\eta)\right),\\
y(\xi,\eta)&=\tfrac{1}{2}+\tfrac{1}{2}\left(\eta+\tfrac{1}{5}\sin(\pi\xi)\sin(\pi\eta)\right).
\end{align}
\end{subequations}

\begin{figure}[htb]
\centering
\includegraphics[width=1.\textwidth]{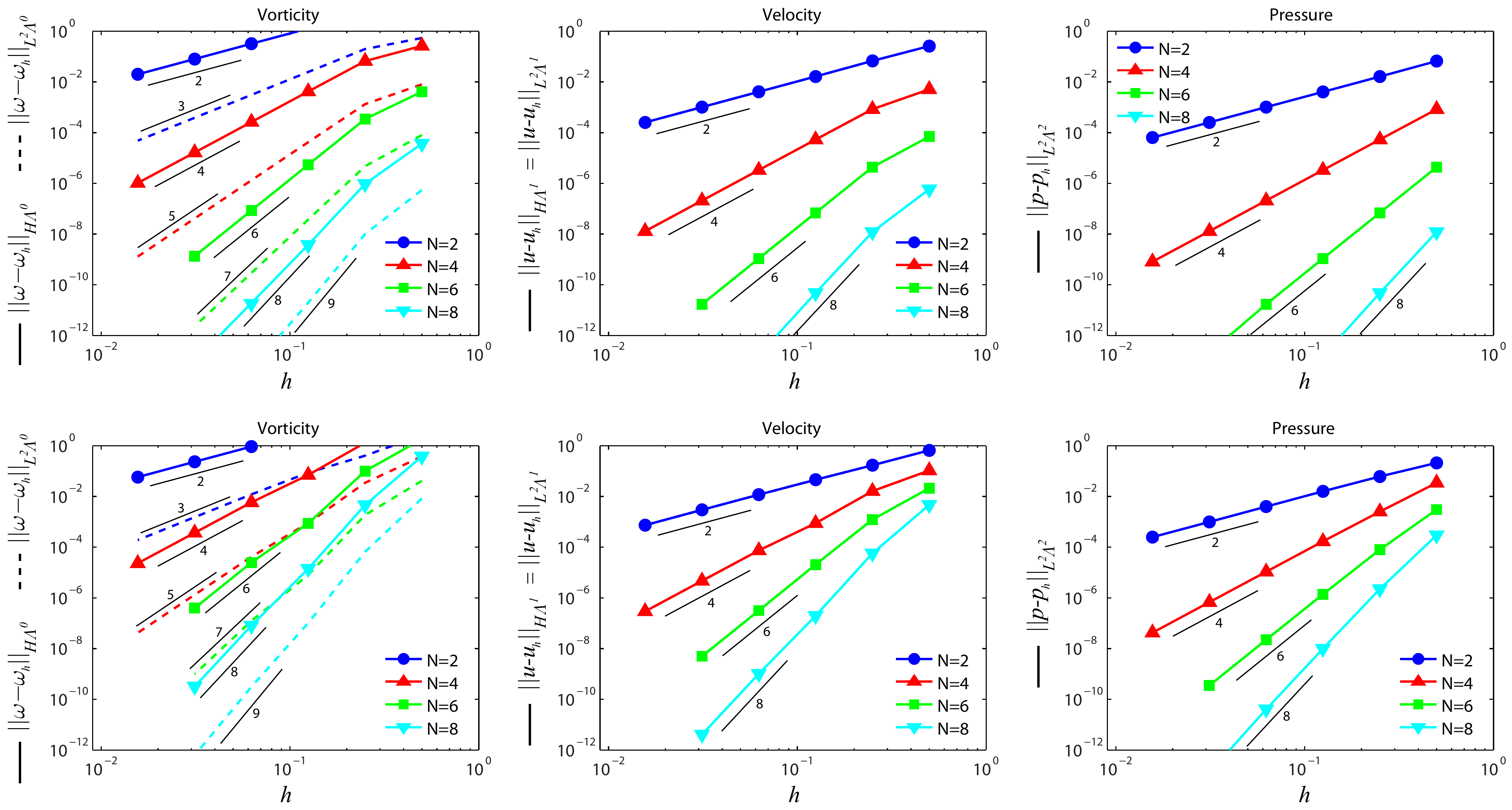}
\caption{Vorticity, velocity and pressure $h$-convergence results of problem \eqref{testcase1}. Results in the top row correspond to Cartesian meshes, results in the bottom row are obtained on curvilinear meshes. All variables are tested on meshes with $N=2,4,6$ and 8.}
\label{fig:hconv}
\end{figure}

\begin{figure}[htb]
\centering
\includegraphics[width=1.\textwidth]{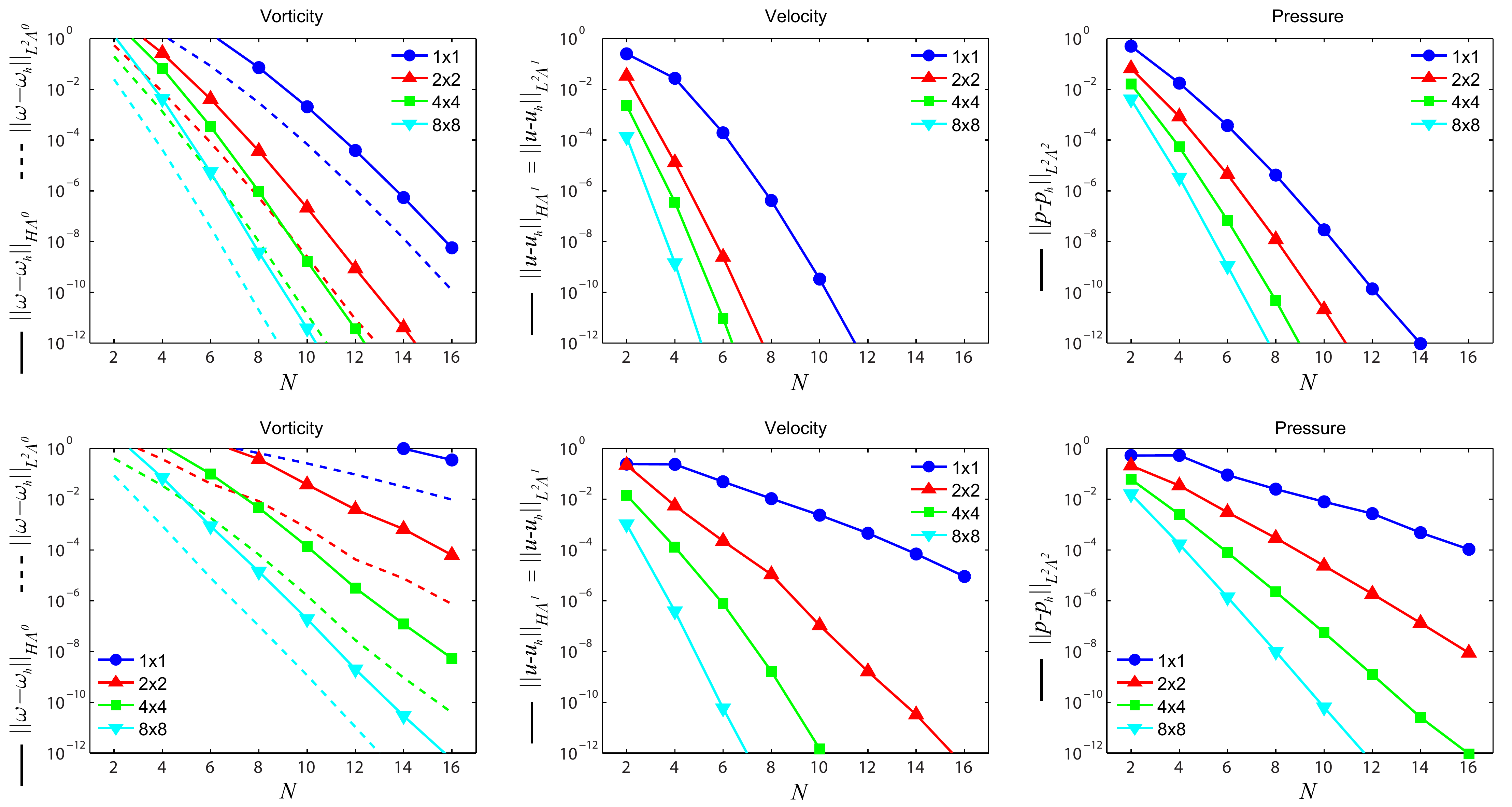}
\caption{Vorticity, velocity and pressure $p$-convergence results of problem \eqref{testcase1}. Results in the top row correspond to Cartesian meshes, results in the bottom row are obtained on curvilinear meshes. All variables are tested on meshes with $1\times1,\ 2\times 2,\ 4\times 4$ and $8\times 8$ spectral elements.}
\label{fig:pconv}
\end{figure}

\begin{figure}[htb]
\centering
\includegraphics[width=0.5\textwidth]{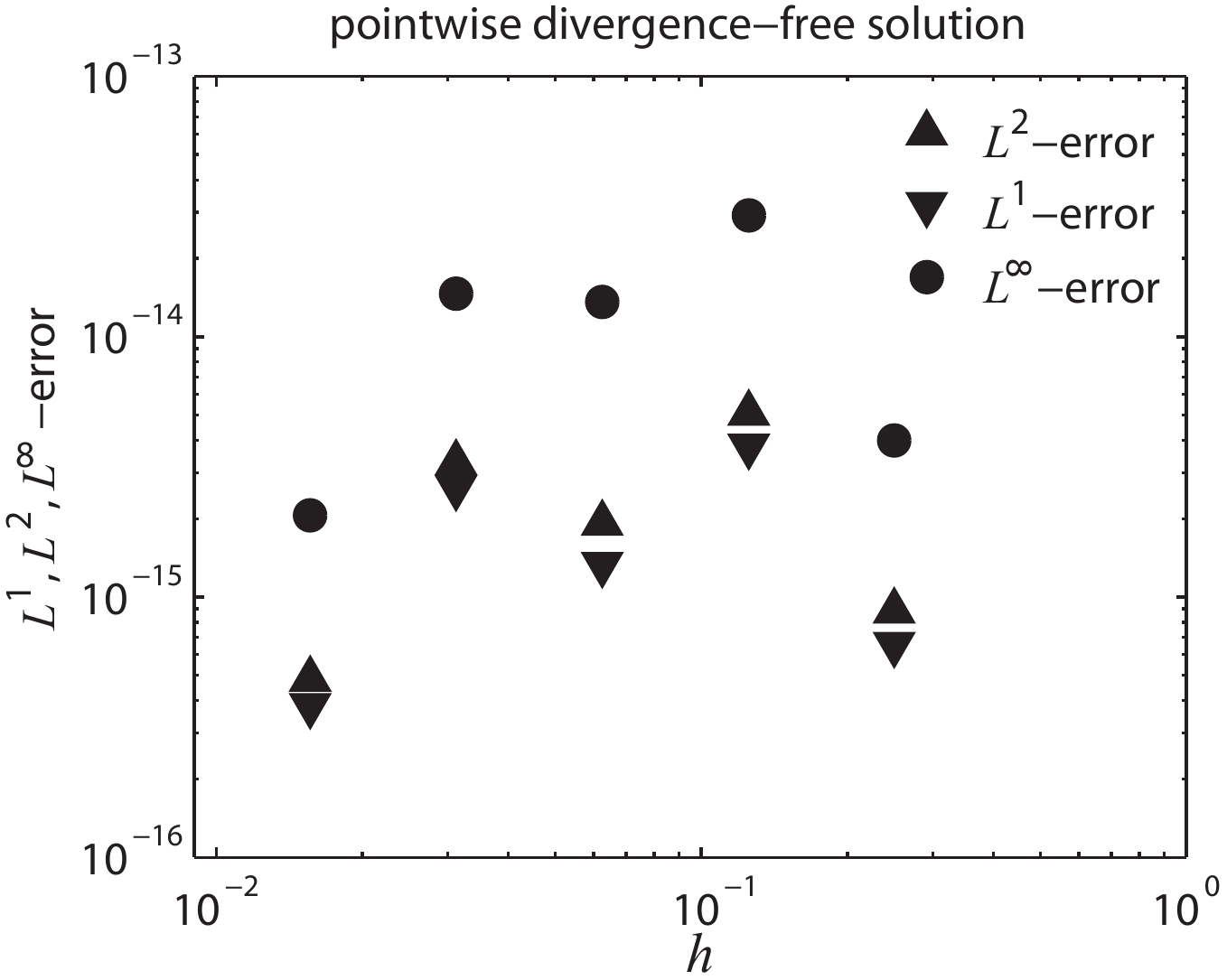}
\caption{$L^1$, $L^2$ and $L^\infty$-error of $\mathrm{div}\,u$ on the Cartesian mesh for discontinuous piecewise linear functions, $N=2$.}
\label{fig:divergencefree}
\end{figure}

\figref{fig:hconv} shows the $h$-convergence and \figref{fig:pconv} shows the $p$-convergence of the vorticity $\kdifformh{\omega}{0}\in\Lambda_h^{0}(\mathcal{Q};C_0)$, velocity $\kdifform{u}{1}_h\in\Lambda^{1}_h(\mathcal{Q};C_1)$ and pressure $\kdifform{p}{2}_h\in\Lambda^2_h(\mathcal{Q};C_2)$. For both figures, the results of the top row are obtained on Cartesian meshes and the results depicted underneath are obtained on curvilinear meshes. The errors for the vorticity and velocity are both measured in the $L^2\Lambda^k$- and $H\Lambda^k$-norm, i.e. $\Vert\omega-\omega_h\Vert_{L^2\Lambda^{0}}$, $\Vert\omega-\omega_h\Vert_{H\Lambda^{0}}$, and $\Vert u-u_h\Vert_{L^2\Lambda^{1}}$, $\Vert u-u_h\Vert_{H\Lambda^{1}}$, respectively. Because the divergence-free constraint is satisfied pointwise, the norm $\Vert\ud( u-u_h)\Vert_{L^2\Lambda^{2}}$ is zero or machine precision, see \figref{fig:divergencefree}, and so the $H\Lambda^{1}$-norm is equal to the $L^2\Lambda^{1}$-norm of the velocity, i.e., $\Vert u-u_h\Vert_{H\Lambda^{1}}=\Vert u-u_h\Vert_{L^2\Lambda^{1}}$. This does not hold for the vorticity, since $\ud\kdifform{\omega}{0}\in\Lambda^{1}_h(\mathcal{Q};C_1)$ is again a function of sine and cosine functions. The norm $\Vert\ud(\omega-\omega_h)\Vert_{L^2\Lambda^{1}}$ converges one order slower than $\Vert\omega-\omega_h\Vert_{L^2\Lambda^{0}}$. More details on the convergence behavior can be found in \cite{kreefterrorestimate}.

In \figref{fig:hconv} the slope of the theoretical convergence rates, \cite{kreefterrorestimate}, are added which shows that $h$-convergence rates are equal to the $h$-convergence rates of the interpolation error \eqref{vvpinterpolationerror}, on both Cartesian as well as curvilinear meshes. \figref{fig:pconv} shows that exponential convergence rates are obtained on both types of meshes.

It is important to remark is that these results are independent of the kind of boundary conditions used. This is shown in Table~\ref{tab:bcresults}. This is an important result, because especially optimal convergence for the normal velocity - tangential velocity boundary condition is non-trivial in compatible methods, \cite{arnold2011}. The standard elements in compatible methods, the Raviart-Thomas elements, show only sub-optimal convergence for velocity boundary conditions, \cite{arnold2011}.

\begin{table}[htb]
\centering
\begin{tabular}{c|c|c|c||c}
 normal velocity & tangential velocity & vorticity & \ vorticity\  & convergence\\
 tangential velocity & pressure & normal velocity & pressure & rate\\ \hline\hline
   4.0758e-01 & 5.4293e-01 & 5.4292e-01 & 5.4292e-01 & $-$  \\
   1.9814e-01 & 1.9738e-01 & 1.9738e-01 & 1.9738e-01 & 1.46 \\
   2.4893e-02 & 2.4776e-02 & 2.4776e-02 & 2.4776e-02 & 2.99 \\
   3.1037e-03 & 3.0954e-03 & 3.0954e-03 & 3.0954e-03 & 3.00 \\
   3.8738e-04 & 3.8684e-04 & 3.8684e-04 & 3.8684e-04 & 3.00 \\
   4.8386e-05 & 4.8352e-05 & 4.8352e-05 & 4.8351e-05 & 3.00 \\
\end{tabular}
\caption{This table shows the vorticity error $\norm{\omega-\omega_h}_{L^2\Lambda^{0}}$ obtained using the four types of boundary conditions described in Table~\ref{tab:boundaryconditions}. The results are obtained on a Cartesian mesh with $N=2$ and $h=\tfrac{1}{2},\tfrac{1}{4},\tfrac{1}{8},\tfrac{1}{16},\tfrac{1}{32},\tfrac{1}{64}$. All four cases show third order convergence.}
\label{tab:bcresults}
\end{table}

\subsection{Lid-driven cavity Stokes}
For many years, the lid-driven cavity flow was considered as one of the classical benchmark cases for the assessment of numerical methods and the verification of incompressible (Navier)-Stokes codes. The lid-driven cavity test case deals with a flow in a unit-square box with three solid boundaries and moving lid as the top boundary, moving with constant velocity equal to one to the right. Because of the discontinuities of the velocity in the two upper corners, the solution becomes singular at these corners, where both vorticity and pressure become infinite. Especially these singularities make the lid-driven cavity problem a challenging test case.

\begin{figure}[htb]
\centering
\includegraphics[width=1.\textwidth]{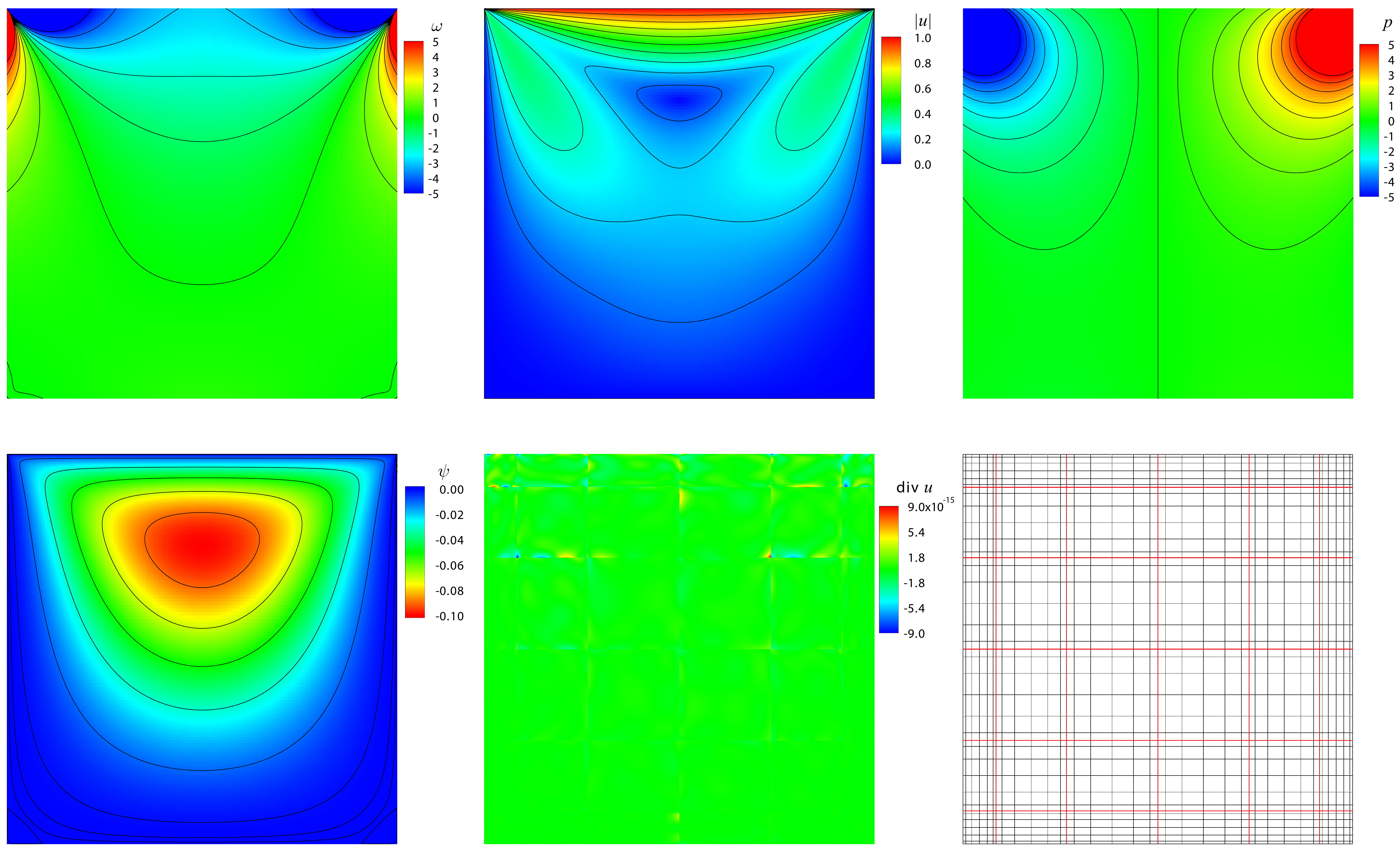}
\caption{Lid-driven cavity Stokes problem results. The top row from left to right shows the solution of the vorticity, velocity magnitude and pressure fields. The bottom row shows from left to right the solution of the stream function, the divergence of the velocity field and the $6\times 6$, $N=6$ mesh.} 
\label{fig:LDC}
\end{figure}

For this test case a non-uniform $6\times6$ Cartesian spectral element mesh is used. Each spectral element consists of a Gauss-Lobatto mesh for $N=6$, see \figref{fig:LDC}. The solutions of the vorticity, velocity, pressure and stream function are shown in \figref{fig:LDC}. Also shown in \figref{fig:LDC} is a plot of the divergence of velocity. It confirms a pointwise divergence-free solution up to machine precision. The results are in perfect agreement with those in \cite{sahinowens2003}.

Because in the mixed mimetic spectral element method no velocity unknowns are located at the upper corners -- only velocity flux {\em through} edges is considered --, no special treatment is needed for the corner singularities, in contrast to many nodal finite-difference, finite-element and spectral element methods, \cite{botellapeyret1998,evans2011,peyrettaylor,prootthesis}. This is due to the finite-volume like structure of the method, as explained in the section of algebraic topology.

\begin{figure}[htbp]
\centering
\includegraphics[width=1.\textwidth]{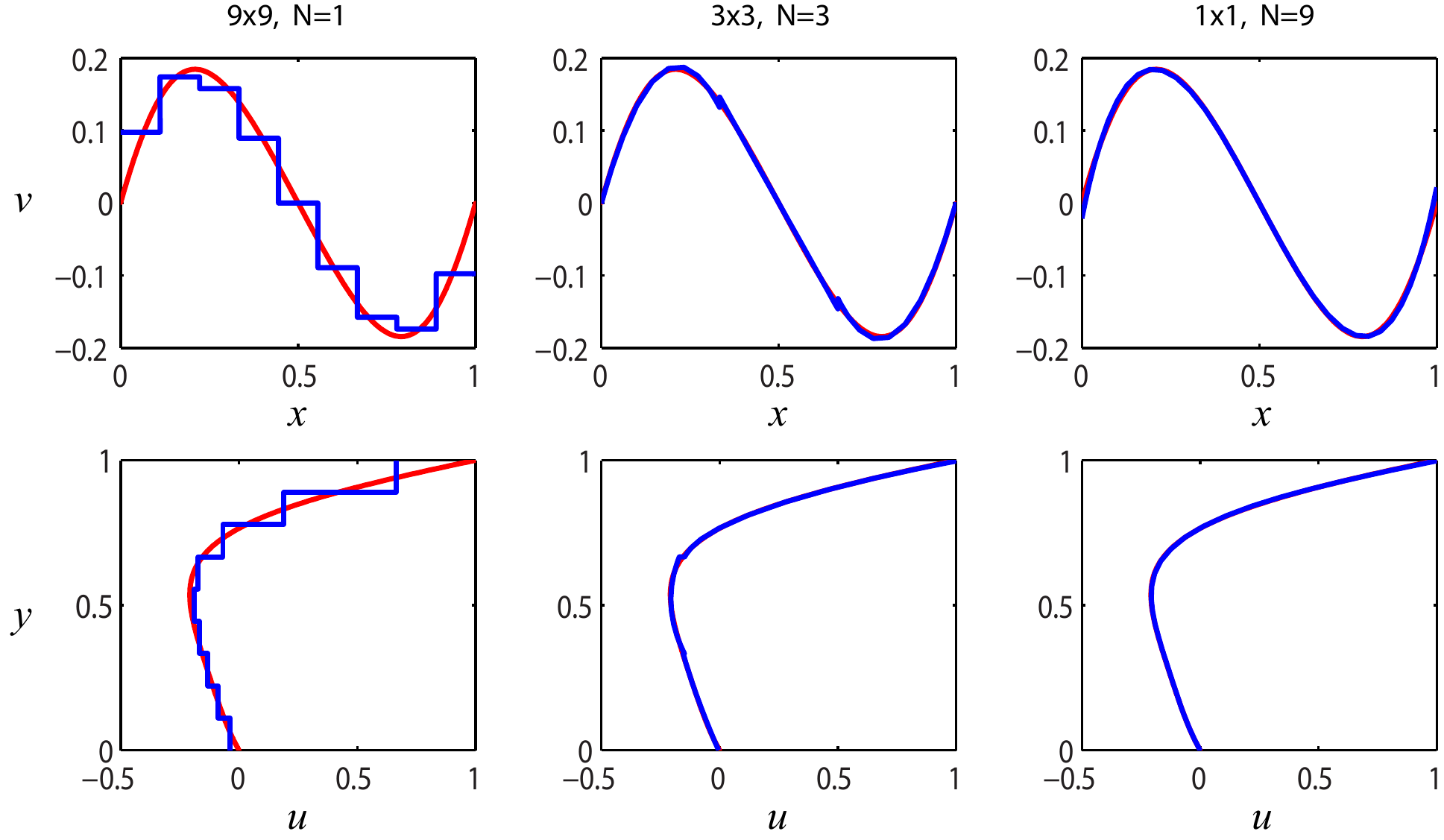}
\caption{Horizontal (top) and vertical (bottom) centerline velocities are shown in blue for a very course mesh, $9\times9$ 2-cells. From left to right the $9\times9$ 2-cells are used in: $9\times9$ zeroth-order elements, $3\times3$ second-order elements and one eight-order element. In red the reference solution from \cite{sahinowens2003}.}
\label{fig:centerlines}
\end{figure}

In \figref{fig:centerlines} the centerline velocities are plotted. Three different configurations are used, based on the same cell complex consisting of $9\times 9$ 2-cells:
\begin{itemize}
\item left: $9\times9$ spectral elements with $N=1$, resulting in piecewise constant approximations along the centerlines,
\item middle: $3\times3$ spectral elements with $N=3$, resulting in piecewise quadratic approximations along the centerlines,
\item right: One global spectral element with $N=9$, resulting in $8^{\rm th}$ order polynomial approximations along the centerlines.
\end{itemize}
Despite the low resolution, all approximations lay on top of those in \cite{sahinowens2003}.

\begin{figure}[htbp]
\centering
\includegraphics[width=1.\textwidth]{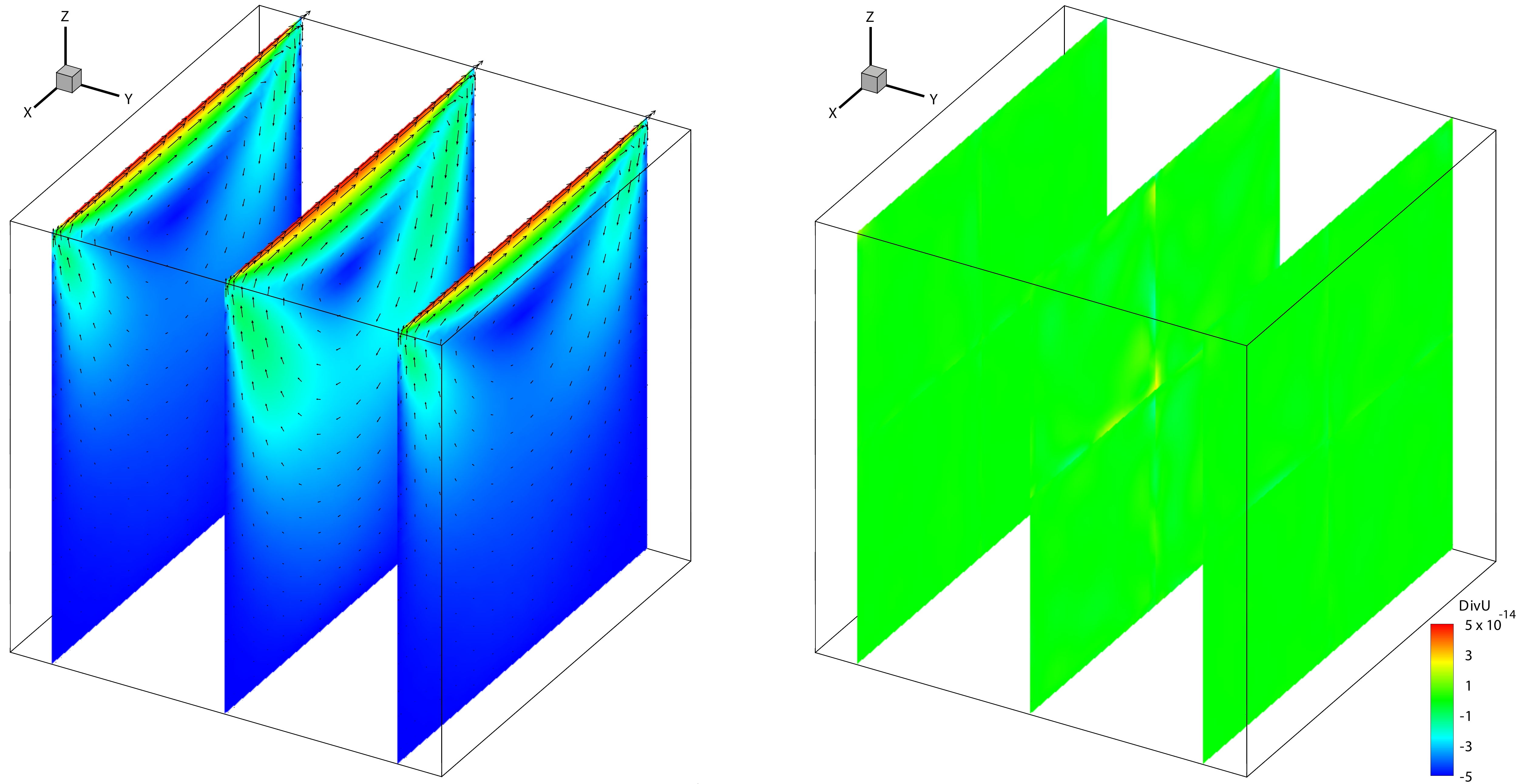}
\caption{Left: slices of magnitude of the velocity field of a three dimensional lid-driven cavity Stokes problem obtained on a $2\times2\times2$ element mesh with $N=8$. Right: slices of the  divergence of velocity. Is confirms a divergence-free velocity field.}
\label{fig:3Dldcdivfree3D}
\end{figure}

Because of the tensor-product construction of discrete unknowns and basis-functions, an extension to three dimensions is straightforward. A 3D lid-driven cavity is of interest because it not only contains corner singularities, but also line singularities. The left plot in \figref{fig:3Dldcdivfree3D} shows slices of the magnitude of the velocity field in a three dimensional lid-driven cavity Stokes problem, obtained on a $2\times2\times2$ element mesh with $N=8$. The slices are taken at 10\%, 50\% and 90\% of the y-axis. The right plot in \figref{fig:3Dldcdivfree3D} shows slices of divergence of the velocity field. The solution at the symmetry plane coincides with the 2D results in \figref{fig:LDC}. It confirms that also in three dimensions the mixed mimetic spectral element method leads to an accurate result with a divergence-free solution.


\begin{figure}[htbp]
\centering
\includegraphics[width=1.0\textwidth]{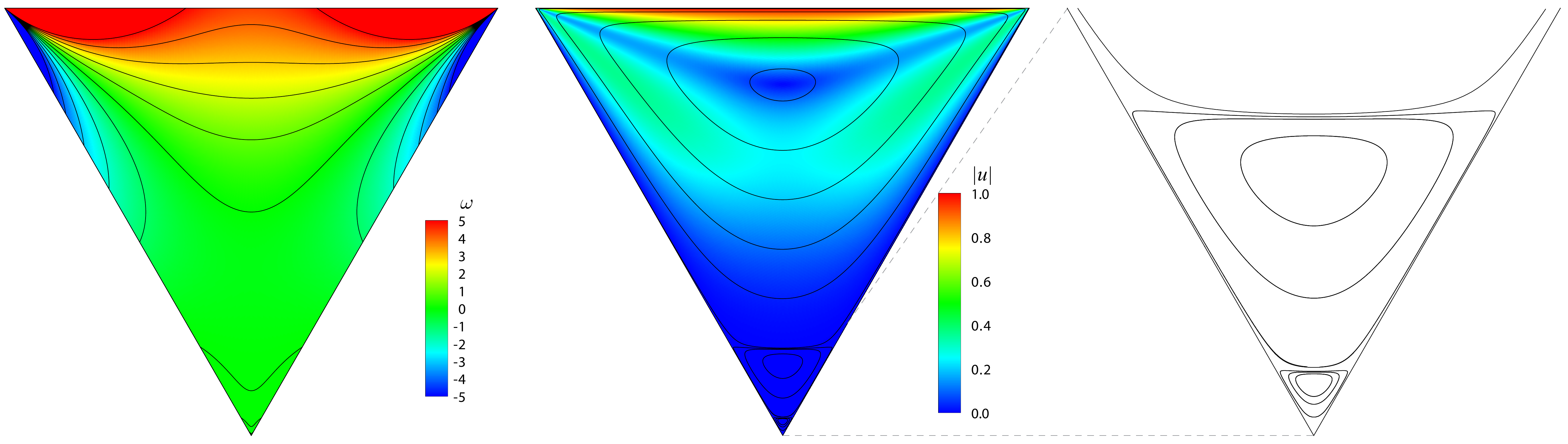}
\caption{Lid-driven cavity Stokes flow in a triangle. Left the vorticity field, in the middle the velocity magnitude with stream function contours on top, and right the stream function contours of a close-up of the bottom corner, revealing the second and third Moffatt eddies.}
\label{fig:triangleLDC}
\end{figure}

The corner singularities can be made even more severe by sharpening the corners, as happens for a lid-driven cavity problem in a triangle. \figref{fig:triangleLDC} shows the vorticity field and the velocity magnitude. On top of the velocity plot, stream function contours are plotted. The solutions are constructed on a 9 spectral element mesh with $N=9$. A close-up of the stream function contours is shown in the rightmost plot in \figref{fig:triangleLDC}. The stream function contours nicely show the first three Moffatt eddies \cite{moffatt}.

\subsection{Flow over a cylinder}
The last test case considers the flow around a cylinder moving with constant velocity to the left, as defined in \cite{changnelson1997}. This testcase is mostly considered in the context of least-squares finite and spectral element methods, due to their moderate performance in case of large contraction regions, \cite{changnelson1997,deanggunzburger,prootgerritsma2006}, mainly in terms of conservation of mass.

\begin{figure}[htbp]
\centering
\includegraphics[width=1.\textwidth]{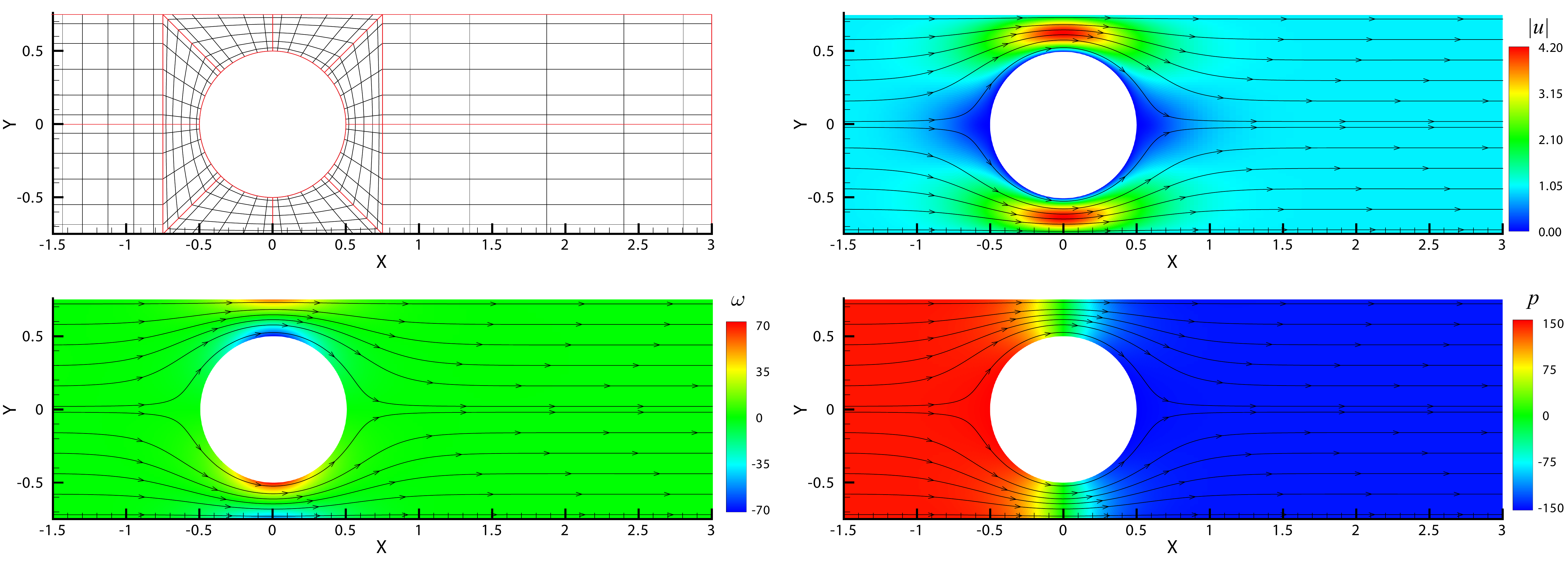}
\caption{Spectral element mesh (top left), magnitude of velocity (top right), vorticity (bottom left) and pressure (bottom right) for flow around a moving cylinder, on a 12 element, $N=6$ mesh.}
\label{fig:cylindercase}
\end{figure}

The cylinder moves with unit velocity along the centerline of a narrow channel. The computational domain is defined as a rectangular box minus the cylinder, as shown in \figref{fig:cylindercase}. Also visible in this figure are the 12 spectral elements in which the computational domain is divided. A transfinite mapping, \cite{gordonhall1973}, is used to define the curved elements around the cylinder. Velocity boundary conditions of $(u,v)=(1,0)$ are prescribed on the outer boundary and no-slip, $(u,v)=(0,0)$, is prescribed along the boundary of the cylinder. Solution of the vorticity, velocity magnitude and pressure, together with streamlines are shown in \figref{fig:cylindercase}.

\begin{figure}[htbp]
\centering
\includegraphics[width=0.4\textwidth]{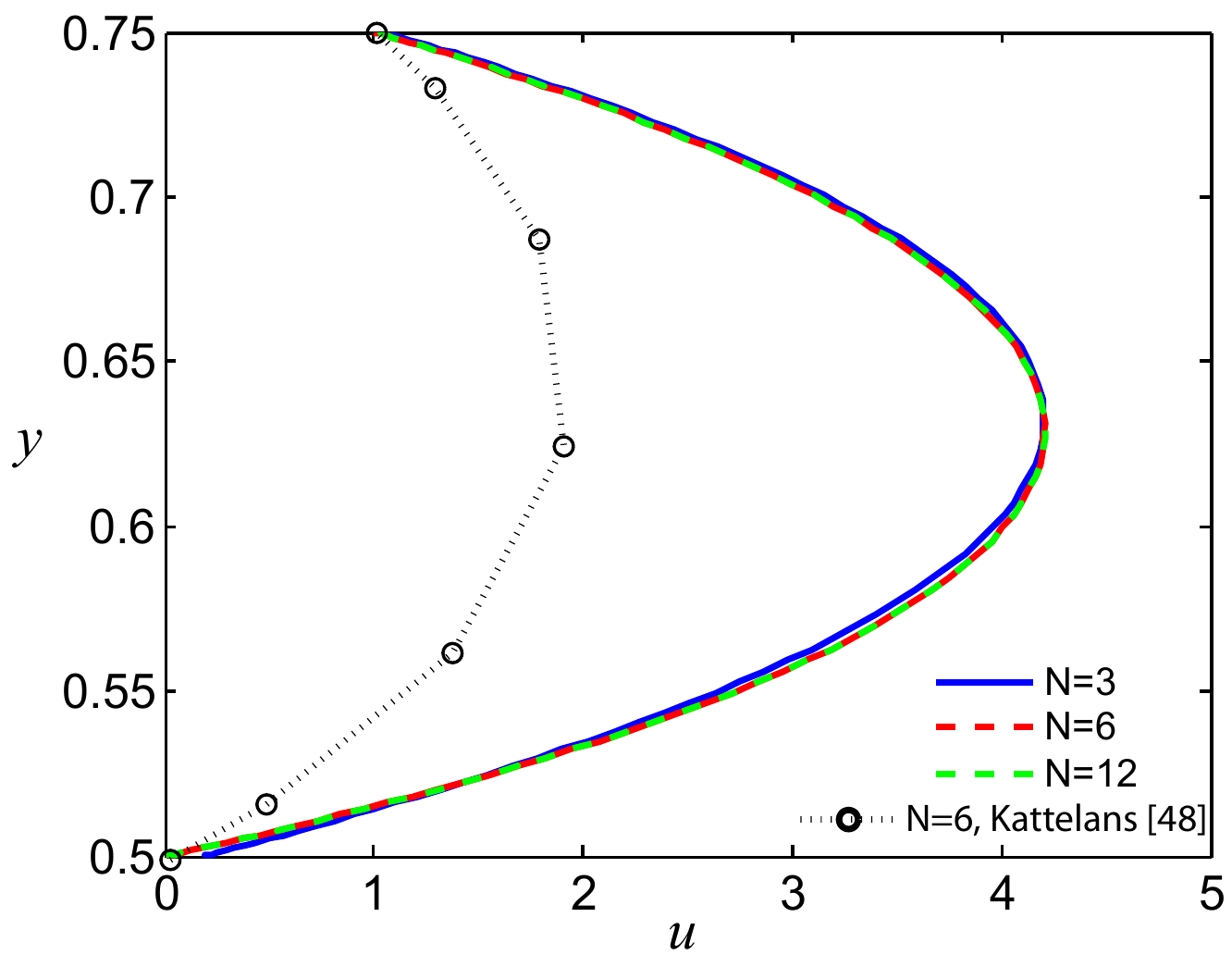}
\caption{Horizontal velocity at smallest cross-section above the cylinder, on a 12 element mesh, for $N=3,6,12$.}
\label{fig:velocity_cross-section}
\end{figure}

Next consider a control volume $\Omega_c$ consisting of the 6 elements in the domain $-1.5\leq x\leq 0,\ 0.75\leq y \leq 0.75$. The control volume is chosen such that the ratio in size between inflow and outflow boundary is maximal. In this control volume conservation of mass should hold. Conservation of mass is expressed, by means of generalized Stokes theorem \eqref{stokestheorem}, in terms of a boundary integral as
\begin{equation}
0=\int_{\Omega_c}\ud \kdifformh{u}{1}\stackrel{\eqref{stokestheorem}}{=}\int_{\p\Omega_c}\kdifformh{u}{1}.
\end{equation}
From Section~\ref{pointwisedivergencefree} and the results of the previous test cases we know that the solution of the velocity is divergence-free throughout the domain, independent of the chosen control volume. In \figref{fig:velocity_cross-section} a comparison is made for the horizontal velocity component $u$ at the smallest cross-section above the cylinder, i.e. $x=0$, $0.5\leq y\leq0.75$, between the recently developed LSSCM, \cite{kattelans2009}, and our MMSEM method for $N=3,6,12$. Both methods use a similar mesh of 12 spectral elements. As can be seen from this figure, the MMSEM method performs already very well for $N=3$, i.e. quadratic polynomial, where the LSSCM still fails for $N=6$, i.e. sixth order polynomial. This is a direct consequence of the pointwise divergence-free discretization.

\section{Conclusions and future aspects}
In this paper we presented the mixed mimetic spectral element method, applied to the vorticity-velocity-pressure formulation of Stokes model. At the heart lies the generalized Stokes theorem, which relates the boundary operator applied on an oriented geometric objects to the exterior derivative, resembling the vector operators grad, curl and div, and the recently developed higher-order mimetic discretization for quadrilaterals and hexadrals, \cite{kreeftpalhagerritsma2011}. The gradient, curl and divergence conforming method results in a point-wise divergence-free discretization of the Stokes problem, as was confirmed by a set of benchmark problems. These results also showed optimal convergence, independent of the type of boundary conditions on orthogonal and curved meshes. More on convergence behavior and error estimates is presented \cite{kreefterrorestimate}. In the near future we plan to extend the method with structure-preserving $hp$-refinement based on a compatible mortar element method.

\bibliographystyle{abbrv}
\bibliography{./literature}

\end{document}